%% file: Klein_Schuster_Wald_arXiv.tex
\documentclass[final,leqno]{siamltex}

\usepackage{graphicx}        
\usepackage{multicol}        
\usepackage[bottom]{footmisc}

\usepackage{newtxmath}       

\usepackage{color}
\usepackage{float}

\newtheorem{algorithm}{Algorithm}
\newtheorem{remark}{Remark}

\author{Rebecca Klein\thanks{Department of Mathematics, Saarland University, PO Box 15 11 50, 66123 Saarbr\"ucken, Germany ({\tt klein@num.uni-sb.de})} \and Thomas Schuster\thanks{Department of Mathematics, Saarland University, PO Box 15 11 50, 66123 Saarbr\"ucken, Germany ({\tt thomas.schuster@num.uni-sb.de})} \and Anne Wald\thanks{Department of Mathematics, Saarland University, PO Box 15 11 50, 66123 Saarbr\"ucken, Germany ({\tt anne.wald@num.uni-sb.de})}}

\title{Sequential subspace optimization for recovering stored energy functions in hyperelastic materials from time-dependent data}

\begin{document}

\maketitle

\begin{abstract}
Monitoring structures of elastic materials for defect detection by means of ultrasound waves (Structural Health Monitoring, SHM) demands for an efficient computation of parameters which characterize their mechanical behavior. Hyperelasticity describes a nonlinear elastic behavior where the second Piola-Kirchhoff stress tensor is given as a derivative of a scalar function representing the stored (strain) energy. Since the stored energy encodes all mechanical properties of the underlying material, the inverse problem of computing this energy from measurements of the displacement field is very important regarding SHM. The mathematical model is represented by a high-dimensional parameter identification problem for a nonlinear, hyperbolic system with given initial and boundary values. Iterative methods for solving this problem, such as the Landweber iteration, are very time-consuming. The reason is the fact that such methods demand for several numerical solutions of the hyperbolic system in each iteration step. In this contribution we present an iterative method based on sequential subspace optimization (SESOP) which in general uses more than only one search direction per iteration and explicitly determines the step size. This leads to a significant acceleration compared to the Landweber method, even with only one search direction and an optimized step size. This is demonstrated by means of several numerical tests. 
\end{abstract}


\input{intro_arXiv.tex}

\input{hyp_mat_arXiv.tex}

\input{sesop_arXiv.tex}

\input{num_res_arXiv.tex}

\input{conclusion_arXiv.tex}

\input{referenc_arXiv}

\end{document}

%% file: intro_arXiv.tex
\section{Introduction}
\label{sec:introduction}

Monitoring structures consisting of materials like fi\-ber-re\-in\-forced plastics or metal laminates is of utmost importance regarding the early detection of defects
such as cracks and delaminations or to estimate the structure's lifetime. Such materials play an important role in the construction of wind power stations, aircrafts and automobiles.
A Structural Health Monitoring (SHM) system consists of a number of actuators and sensors that are applied to the structure. We refer to the seminal book of Giurgiutiu \cite{Giurgiutiu2008} 
for a comprehensive outline of piezoelectric sensor based SHM systems and their mechanics. A comprehensive monograph on Lamb wave based SHM in polymer composites
is given by \cite{LAMBWAVES-BUCH:18}. The mechanical waves that are generated
by the actuators propagate through the structure, interact with a possible damage and are measured at the sensors. The inverse problem then consists in recovering
the damage from the given sensor measurements. The mathematical model of wave propagation in solids is represented by Cauchy's equation of motion
\[  \rho \ddot{u} - \nabla\cdot P = f,   \]
where $\rho$ denotes the mass density, $P$ the first Piola-Kirchhoff stress tensor, $f$ an external volume force vector and $u$ is the displacement field of the wave.
Materials such as fiber-reinforced plastics or metal laminates are elastic and, depending on the respective response function for $P$, we obtain a corresponding system of
hyperbolic partial differential equations for the displacement field $u$. The response function for $P$ in turn encodes macroscopic mechanical properties of
the material, such as, e.g., the Poisson number or Young's modulus, yielding pointers to hidden damages. There is a vast amount of literature concerning inverse problems connected
to Cauchy's equation of motion in elasticity and we refer here only to recent works that have a close relation to the topic of this contribution. Inverse problems
in linear elasticity are, e.g., considered in \cite{BINDER;ET;AL:15, BOURGEOIS;ET;AL:11, HUBMER;ET;AL:18, LECHLEITER;SCHLASCHE:17}. A nice overview on inverse problems in
elasticity is \cite{BONNET;CONSTANT:05}. 
An important material class in elasticity is given by \emph{hyperelastic} materials, which are characterized
by the fact that the stress tensor is given as a derivative of a scalar function with respect to the strain tensor. This scalar function is the stored (strain) energy
function and its integral equals the total strain energy which is necessary to deform the body. Since all relevant material properties can be deduced from the stored energy
function, its computation should reveal valuable pointers to damages in the structure. The corresponding Cauchy equation then is nonlinear. In \cite{RAUTER;LAMMERING:15}
the authors investigate higher harmonics of Lamb waves in hyperelastic isotropic materials. Inverse problems in nonlinear elasticity are, e.g., considered in
\cite{DEHOOP;ET;AL:18, SCHUSTER;WOESTEHOFF:14, SEYDEL;SCHUSTER:16, SEYDEL;SCHUSTER:17, SRIDHAR;ET;AL:18}. In the present contribution we consider the nonlinear inverse problem
of reconstructing the stored energy function from the knowledge of the full displacement field $u$. The stable solution of nonlinear, dynamic inverse problems is currently 
counted among the most demanding mathematical challenges.

Nonlinear inverse problems are usually solved by iterative regularization techniques. Standard methods such as the Landweber iteration scheme prove to be tremendously slow
when applied to such a high-dimensional nonlinear inverse problem. To increase numerical efficiency Sequential Subspace Optimization (SESOP) techniques have been developed and 
analyzed for various settings, see \cite{Narkiss, SchopferSchuster, schopfer2008metric, WaldSchuster, Wald2018}. The general idea is to reduce the number of iterations until the stopping criterion is fulfilled. To this end, the classical Landweber method is extended by two features. First, a finite number of search directions is used in each iteration. Second, the length of each search direction is explicitly calculated. This is done in such a way that the method admits a very intuitive interpretation: The iterate is sequentially projected onto subsets that contain the solution set of the inverse problem. These subsets are intersections of stripes that correspond to the respective search directions. The calculation of the projection yields a regulation of the step widths. \\
This technique has been successfully applied, for example in parameter identification \cite{schuster2012regularization, awts18, Wald2018}, demonstrating a significant increase in efficiency. In addition, they have been used and analyzed in combination with, e.g., sparsity constraints, total variation or Nesterov methods \cite{GuHanChen19, MaassStrehlow16, TongHanLongGu19}.

This contribution delivers a proof-of-concept by demonstrating that RESESOP applied to a high-dimensional nonlinear and dynamic inverse problem 
leads to a significantly faster convergence as well as less computation time with at the same time higher accuracy compared to Landweber's method. 

\emph{Outline.} In Section \ref{sec:hyperelasticMaterials} we briefly summarize essential concepts of continuum mechanics for elastic solids and deduce the exact mathematical
setting for identifying the stored energy function of a hyperelastic material from measurements of the displacement field. In order to guarantee that the reconstructed
energy is physically meaningful we use a dictionary of finitely many elements. The inverse problem subsequently reduces to the computation of the corresponding
coefficients with respect to the given dictionary. Section \ref{sec:results} outlines the introduction and analysis of the Landweber method and RESESOP. In Section
\ref{sec:num_res} we finally present several numerical experiments using three different damage scenarios for a structure consisting of a Neo-Hookean material
showing the superiority of RESESOP compared to the Landweber method.

%% file: hyp_mat_arXiv.tex
\section{Hyperelastic materials}
\label{sec:hyperelasticMaterials}

In this chapter we briefly discuss some basic facts from continuum mechanics and especially on Cauchy's equation of motion and hyperelastic constitutive equations.
For deeper insights we refer to the standard literature \cite{Ciarlet, holzapfel2000nonlinear, marsden1994mathematical}.

The considered elastic structure is described by a bounded, open, connected subset $\Omega\subset \mathbb{R}^3$ with a sufficiently smooth boundary. We
start with the mathematical definition of a deformation of $\Omega$.\\

\begin{definition}\label{dfn:deformation}
A \emph{deformation} of a body $\Omega$ is an invertible, continuously differentiable mapping $\varphi:[0,T]\times \Omega \to \mathbb{R}^3$, which is orientation-preserving such that
\begin{equation*}
\det (\nabla \varphi(t,x))>0 \qquad \forall (t,x) \in [0,T] \times \Omega,
\end{equation*}
where
\begin{equation*}
\nabla \varphi (t,x) = \left( \frac{\partial \varphi_i}{\partial x_j}(t,x) \right)_{i,j=1,2,3} = \left( \partial_{x_j} \varphi_i (t,x)\right)_{i,j=1,2,3} \in \mathbb{R}^{3\times 3}.
\end{equation*}
\end{definition}

Definition \ref{dfn:deformation} implies that the body will not be torn apart or penetrate itself during the deformation. Since $\nabla \varphi$ is invertible, any two points in $\Omega$
can be separated at any time $t\in [0,T]$. If $\Omega$ undergoes a deformation, then a fixed point $x$ is shifted to a point $\varphi (t,x)$. Their difference defines
the displacement field.\\

\begin{definition}
Let $\varphi:[0,T] \times \Omega \to \mathbb{R}^3$ be a deformation.  
Then the \emph{displacement field} $u:[0,T] \times \Omega \to \mathbb{R}^3$ is given by 
\begin{equation*}
u(t,x)=\varphi(t,x)-x .
\end{equation*} 

\end{definition}

The set $\Omega$ is also called the \emph{reference configuration} whereas $\Omega (t):= \varphi (t,\Omega)\subset \mathbb{R}^3$ is called the \emph{deformed configuration} and represents the body
after deformation at time $t$. A guided wave that is generated by actuators and propagates through the structure $\Omega$ will cause a displacement field $u$ which subsequently can be measured by applied sensors. This is the key idea of an SHM system (c.f. \cite{Giurgiutiu2008}).\\

\begin{definition} \label{def:verschiebungsgradient}
Let $\varphi:[0,T] \times \Omega \to \mathbb{R}^3$ be a deformation and $u$ be the corresponding displacement field. 
The \emph{displacement gradient} is given by
\begin{equation*}
\nabla u(t,x)= \nabla \varphi(t,x) - I
\end{equation*}
with the identity matrix $I\in \mathbb{R}^{3\times 3}$. The gradient $\nabla$ refers to the spatial coordinates.\\
\end{definition}

The propagation of ultrasound waves in $\Omega$ is mathematically described by \emph{Cauchy's equation of motion}, which follows from
the stress principle of Euler and Cauchy and the axioms of force and moment balance.
For all $(t,x) \in [0,T] \times \Omega$ we have 
\begin{equation} \label{dfg first version}
\rho(x) \ddot{u}(t,x) -\nabla \cdot P(t,x) =f(t,x).
\end{equation}
Here $\rho:\Omega \to \mathbb{R}^+$ denotes the mass density, $f:[0,T] \times \Omega \to \mathbb{R}^3$ the external body force and $P:[0,T]\times \Omega \to \mathbb{R}^{3\times 3}$ the first Piola-Kirchhoff stress tensor. This is a differential equation for the unknowns $u$ and $P$ and obviously not uniquely solvable in its present form. But by now we did not include the phenomenon
of elasticity to $\Omega$ and equation \eqref{dfg first version}. Elasticity means that there is a stress-strain relation which is implied by the existence of a so called \emph{response function} for the Cauchy stress tensor. To be short: a deformation of the body $\Omega$ causes strain which again causes stress. Postulating the existence of a response function will furthermore reduce the degrees of freedom in \eqref{dfg first version}.

Before we formulate the principle of elasticity we introduce by $\sigma:[0,T]\times \Omega(t) \to \mathbb{R}^{3 \times 3}$ the \emph{Cauchy stress tensor}. This is a continuously differentiable,
symmetric tensor field whose existence follows from Cauchy's theorem. In some sense this is the counterpart to the first Piola-Kirchhoff stress tensor $P$: The Cauchy stress tensor $\sigma$ is defined
on the deformed configuration $\Omega (t)$ whereas $P$ is defined on the reference configuration $\Omega$. Of course, each of these can be transformed into the respective other one, and the specific relation between $\sigma$ and $P$ is given by equation \eqref{Piola_P}.\\


\begin{definition} \label{dfn:elastisch}
A material is called \emph{elastic}, if a mapping 
\begin{equation*}
\tilde{\sigma } : \overline{\Omega} \times \text{GL}_+ (3) \to \text{Sym}(3),\qquad  (x,Y) \mapsto \tilde{\sigma}(x,Y)
\end{equation*}
exists, such that the Cauchy stress tensor satisfies
\begin{equation} \label{glg: sigma tilde}
\sigma(t,\varphi(t,x))= \tilde{\sigma}(x, \nabla \varphi(t,x))
\end{equation} 
for every deformation $\varphi$, where 
\begin{equation*}
\text{GL}_+(3):= \{ Y \in \mathbb{R}^{3 \times 3 } | \det(Y) > 0 \}
\end{equation*} 
denotes the set of $3 \times 3$ matrices with a positive determinant and Sym(3) is the set of symmetric $3\times 3$ matrices.
The function $\tilde{\sigma}$ is called the \emph{response function} for $\sigma$. Equation \eqref{glg: sigma tilde} is called a \emph{constitutive equation of the material}.\\
\end{definition}


The first Piola-Kirchhoff stress tensor can be computed from $\sigma$ by applying the \emph{Piola transform}

\begin{equation}\label{Piola_P} P(t,x)= \det(\nabla\varphi(t,x)) \sigma(t,x) \nabla \varphi(t,x)^{-\top} . \end{equation}

So, if there exists a response function $\tilde{\sigma}$ for $\sigma$, then we easily obtain a response function $\tilde{P}$ for $P$ from \eqref{Piola_P} via

\begin{equation*}
	 \tilde{P}(x,Y) := \det Y \tilde{\sigma}(x,Y) Y^{-\top},\qquad x\in\Omega,\; Y\in \text{GL}_+ (3).
	 \end{equation*}

\begin{remark}
The \emph{Cauchy-Green strain tensor} $B$ is defined as $B=\nabla \varphi^{\top} \nabla\varphi$ and we have $B=I$ if and only if the deformation is rigid. Thus, $B$ measures the 'deviation'
between a deformation $\varphi$ and a rigid motion. It is quite obvious from \eqref{glg: sigma tilde}, that the existence of a response function $\tilde{\sigma}$ implies the existence of a function $\hat{\sigma}$ with
\[   \sigma (t,x) = \tilde{\sigma} (x,\varphi (t,x)) = \hat{\sigma}(x,B(t,x)).   \]
In this way \eqref{glg: sigma tilde} can be interpreted as a relation between stress and strain which is the reason why \eqref{glg: sigma tilde} is also called \emph{stress-strain relation}.
Hence, elasticity in fact means that a material replies to strain with stress.\\
\end{remark}

A large class of physically very important elastic materials is represented by the \emph{hyperelastic materials}. For this class the response functions have a very specific form.\\

\begin{definition} \label{def: hyperel}
An elastic body is called \emph{hyperelastic} if the response function of the first Piola-Kirchhoff stress tensor is given by 
\begin{equation*}
\tilde{P}(x,Y)=\nabla_Y \hat{C}(x,Y),\qquad x\in\Omega,\; Y\in M,
\end{equation*}
and a scalar function $\hat{C}:\Omega \times \text{GL}_+(3) \to \mathbb{R}$. This function $\hat{C}$ is called \emph{stored (strain) energy function}.\\
\end{definition}

The derivative $\nabla_Y$ used in Definition \ref{def: hyperel} is to be understood as
\begin{equation*}
\nabla_Y g(x,Y) = \left[ \partial_{Y_{i,j}} g(x,Y) \right]_{1 \leq i,j \leq 3} \in \mathbb{R}^{3\times 3},\qquad x\in\Omega,\; Y\in \text{GL}_+ (3),
\end{equation*}
for a differentiable function $g:\Omega \times M \to \mathbb{R}$ and $M\subset \mathbb{R}^{3\times 3}$.\\

\begin{remark}
a) If $\varphi$ is a deformation and the body $\Omega$ consists of a hyperelastic material, then the integral
\[   E(t) = \int_{\Omega} \hat{C}\big( x,\nabla \varphi (t,x) \big)\, \mathrm{d} x   \]
denotes the \emph{strain energy} $E(t)$ which is necessary to perform the deformation at time $t$. This explains the term stored (strain) energy function for $\hat{C}$.\\
b) The fourth order \emph{elasticity tensor} $\mathbb{C}$ can directly be computed from $\hat{C}$ by
\[  \mathbb{C} (x) = \nabla_Y \nabla_Y \hat{C}(x,I),\qquad x\in\Omega .   \]
It plays a crucial role in linear elasticity and its entries are important functions describing material properties such as Young's modulus and the Poisson number.
In this sense $\hat{C}$ encodes all important material properties and yields pointers for defects in hyperelastic structures.\\
\end{remark}

Let $\Omega$ be hyperelastic. Then Cauchy's equation of motion reads
\begin{equation}\label{Cauchy-hyper}
\rho(x)\ddot{u}(t,x) - \nabla \cdot \nabla_Y \hat{C}(x, \nabla u(t,x)) = f(t,x), \qquad (t,x) \in [0,T]\times \Omega.
\end{equation}

Note that in \eqref{Cauchy-hyper} we silently used the identity $\nabla u = \nabla\varphi - I$ to write, in slight misuse of notation, $\hat{C}(x, \nabla u(t,x))$.
This means that, by assuming $\Omega$ to be hyperelastic and $\hat{C}$ to be known explicitly, Cauchy's equation of motion is no longer underdetermined since we have three equations and three unknowns,
i.e., the three components of the displacement vector $u$. To ensure uniqueness one furthermore has to postulate initial and boundary values for $u$ (c.f. \cite{WOESTEHOFF;SCHUSTER:14}).

The inverse problem which is numerically solved in this contribution consists in computing the stored energy function $\hat{C}$ from measurements of the displacement field $u$. 
To specify this we follow the idea of computing $\hat{C}$ as a conical combination with respect to a given dictionary consisting of physically reasonable stored energy functions 
$C_K$, $K=1,\ldots,N$., c.f. \cite{KALTENBACHER;LORENZI:07, SCHUSTER;WOESTEHOFF:14}. Let $\{C_K : \Omega\times \mathbb{R}^{3\times 3} \to \mathbb{R} : K=1,\ldots,N \}$ be such a dictionary. Then we write
\begin{equation*}
\hat{C}(x,Y)= \sum_{K=1}^{N} \alpha_K C_K(x,Y),\qquad x\in\Omega,\ Y\in\mathbb{R}^{3\times 3},
\end{equation*}
for certain coefficients $\alpha_K\geq 0$. Equipped with appropriate initial and boundary values we obtain Cauchy's equation of motion in its final form: The balance equation reads
\begin{equation} \label{system:diff}
\rho \ddot{u}(t,x)- \sum_{K=1}^{N} \alpha_K \nabla \cdot \nabla_Y C_K (x, \nabla u(t,x)) = f(t,x),\qquad (t,x) \in [0,T]\times \Omega.
\end{equation}
We furthermore assume initial values
\begin{eqnarray} \label{system:diff2}
u(0,\cdot)&= u_0 \in H^2(\Omega, \mathbb{R}^3),  \\
\dot{u}(0,\cdot)&=u_1 \in H^1(\Omega, \mathbb{R}^3)
\end{eqnarray}
as well as homogeneous boundary values
\begin{equation} \label{system:diff3}
u(t,\xi)=0, \quad \xi \in \partial\Omega.
\end{equation}

The respective inverse problem is formulated as follows:\\

\textbf{(IP)} Given $(f,u_0,u_1)$ and the displacement field $u(t,x)$ for $t\in [0,T]$ and $x\in \Omega$, determine the coefficients $\alpha=(\alpha_1, ... , \alpha_N)\in \mathbb{R}^n_+$, such that $u$ satisfies the initial boundary value problem \eqref{system:diff}--\eqref{system:diff3}.

If we define by $F: \mathcal{D}(F)\subset \mathbb{R}_+^N \to \mathcal{X}$ the forward operator which maps, for fixed given $(f,u_0,u_1)$, a vector $\alpha \in \mathbb{R}_+^N$ to the unique solution
$u\in \mathcal{X}$, then the inverse problem demands for solving the nonlinear operator equation
\begin{equation*}
F(\alpha)= u.
\end{equation*}
Here $\mathcal{D}(F)$ denotes the domain of $F$ consisting of those $\alpha\in \mathbb{R}_+^N$ admitting a unique solution and 
$\mathcal{X}=L^\infty([0,T]\times \Omega,\mathbb{R}^3)\cap W^{1,\infty}(0,T;H^1 (\Omega,\mathbb{R}^3))$ denotes the image space of $F$ 
containing all admissible solutions. For more details regarding existence and uniqueness of solutions for the IBVP \eqref{system:diff}--\eqref{system:diff3}
we refer the reader to \cite{SEYDEL;SCHUSTER:16, WOESTEHOFF;SCHUSTER:14}.

In Section \ref{sec:num_res} we will see that a convenient approach to define the dictionary elements $C_K$ is to use
tensor products
\begin{equation*}
C_K(x,Y)=v_K(x)\hat{C}(Y),\qquad K=1,\ldots,N,
\end{equation*}
with B-splines $v_K$ that are also used for the Finite Element solution of \eqref{system:diff} and physically reasonable stored energy functions $\hat{C}$ depending only on $Y$.
This idea is taken from \cite{SEYDEL;SCHUSTER:17}.

%% file: sesop_arXiv.tex
\section{Sequential subspace optimization}
\label{sec:results}

In this contribution we present numerical results that are obtained with both the attenuated Landweber as well as the RESESOP method as a solver for the inverse problem (IP). 
In \cite{SEYDEL;SCHUSTER:17} some results using the attenuated Landweber method, implemented in \verb#C++# together with the finite element library \verb#deal.II# \cite{dealII}, have already been presented. For the reader's convenience, we will introduce some notation and briefly summarize the attenuated Landweber method. 

\vspace*{2ex}


Consider a (nonlinear) problem 
\begin{align*} 
 F(x)=y, \quad F: \mathcal{D}(F) \subset X \rightarrow Y,
\end{align*}
with Hilbert spaces $X$ and $Y$. Then the respective attenuated Landweber iteration reads
\begin{align}\label{eq_landweber}
x_{k+1}^\delta = x_k^\delta + \omega F'(x_k^\delta)^*(y^\delta -F(x_k^\delta)),\qquad k=0,1,\ldots
\end{align} 
where the parameter $\omega > 0$ is called a \emph{relaxation} or \emph{damping parameter}. Since $\omega$ is fixed, there is no strategy to adapt the step width in each individual iteration.  
It is assumed that we only have disturbed data $y^\delta$ with $ \lVert y^\delta -y \rVert < \delta $ and noise level $\delta >0 $ at our disposal.  
The convergence of the Landweber method is guaranteed by selecting 
\begin{align*}
\omega \in \left(0, \frac{1}{C_\rho^2} \right)
\end{align*}
with the constant 
\begin{align*}
C_\rho:= \sup \lbrace\lVert F'(x) \rVert : x \in B_\rho(x_0) \rbrace.
\end{align*}   
In case of noisy data, the iteration is stopped by the discrepancy principle, which turns it into a regularization method \cite{KaltenbacherNeubauerScherzer, scherzer98}. 

However, the Landweber method is known to be very slowly converging, and it often takes a lot of iterations to obtain a suitable regularized solution. Particularly in view of an application in parameter identification, where the calculation of each gradient involves the numerical evaluation of the forward operator as well as the adjoint of its linearization, a reconstruction via the Landweber method is too time-consuming and hardly practicable, see, e.g., \cite{SEYDEL;SCHUSTER:17, awts18}. 

In contrast to the attenuated Landweber method, the SESOP method not only involves a regulation of the step width, it also potentially uses multiple search directions per iteration. This of course requires additional (but numerically cheap) calculations in each iteration step, such that a SESOP step will take slightly longer. However, we anticipate that the SESOP and RESESOP methods will need far less iterations and thus lead to a faster convergence of the iteration. 

In this section we will give a short introduction to sequential subspace optimization (SESOP) and regularizing sequential subspace optimization (RESESOP). From the RESESOP method we derive the algorithm which we will use for our later experiments, where we solve (IP) numerically from simulated noisy data.

The idea behind the SESOP method and its regularizing version RESESOP is to reduce the number of iteration steps by sequentially projecting the current iterate onto suitable subsets of the source space $X$ that are hyperplanes or stripes in $X$ and contain the solution set of the respective inverse problem $F(x) = y$. This approach is inspired by the fact that in the case of linear problems, the solution set itself is an affine subspace.
More detailed information about the SESOP method for linear problems can be found in \cite{Narkiss, SchopferSchuster,schopfer2008metric}. Results concerning the SESOP method as a solution technique for nonlinear problems are presented in \cite{GuHanChen19, TongHanLongGu19, WaldSchuster, Wald2018}.

\subsection{Basics}
We will first state some basics for the RESESOP method, in particular the definitions of hyperplanes, half-spaces and stripes, as well as the metric projection.\\

\begin{definition}{Hyperplanes, half-spaces and stripes}\\
Let $u \in X \setminus \lbrace 0 \rbrace$ and $\alpha, \xi \in \mathbb{R}$, $\xi \geq 0$. For these parameters, we define the \emph{hyperplane}
  \begin{align*}
   H(u,\alpha) := \left\lbrace x \in X \ : \ \left\langle u,x \right\rangle = \alpha \right\rbrace,
  \end{align*}
  the \emph{half-space}
 \begin{align*}
   H_{\leq}(u,\alpha) := \left\lbrace x \in X \ : \ \left\langle u,x \right\rangle \leq \alpha \right\rbrace,
  \end{align*}
  and the \emph{stripe}
 \begin{align*}
   H(u,\alpha,\xi) := \left\lbrace x \in X \ : \ \left\lvert \left\langle u,x \right\rangle - \alpha \right\rvert \leq \xi \right\rbrace.
  \end{align*}
\end{definition}

The half-spaces $H_\geq(u,\alpha)$, $H_<(u,\alpha)$ and $H_>(u,\alpha)$ are defined analogously. We see that the half space $H_<(u,\alpha)$ is simply the space beneath the hyperplane $H(u,\alpha)$. 
The stripe $H(u,\alpha,\xi)$ emerges from the hyperplane $H(u,\alpha)$ by admitting a width that is determined by $\xi$. 
Hyperplanes, half-spaces as well as stripes are convex, non-empty sets according to their definition. In addition, the sets $H(u,\alpha)$, $H_{\leq}(u,\alpha)$, $H_{\geq}(u,\alpha)$ and $H(u,\alpha,\xi)$ are closed. \\

The solution set $M_{Fx=y}$ of a linear operator equation $Fx=y$ can be described by
\begin{displaymath}
 M_{Fx=y} := \left\lbrace x \in X \, : \, Fx=y \right\rbrace = x_0 + \mathcal{N}(F)
\end{displaymath}
for some $x_0 \in \mathcal{N}(F)^{\bot}$. 

Another tool that plays an important role is the metric projection.\\

\begin{definition}\label{dfn:metrische Projektion}
The \emph{metric projection} of $x\in X$ onto a non-empty closed convex set $C\subset X$ is the unique element $P_C(x) \in C$, such that
\begin{align*}
\lVert x- P_C(x)\rVert^2 =\min_{z\in C} \lVert x-z \rVert^2.
\end{align*}
\end{definition}

The metric projection $P_C$ onto a convex set fulfills the descent property of the form
\begin{align}\label{eq_descent_property}
\lVert z-P_C(x) \rVert^2 \leq \lVert z-x \rVert^2 - \lVert P_C(x) -x \rVert^2
\end{align}
for all $z\in C$. 


Since hyperplanes and stripes are, by definition, closed and convex non-empty sets, the metric projection of $x \in X$ onto these specific subsets is well-defined. For example, if $C := H(u,\alpha)$ is a hyperplane of $X$, then the metric projection of $x \in X$ onto $C$ corresponds to the orthogonal projection, i.e., we have
\begin{align} \label{proj}
P_{H(u,\alpha)}(x)=x-\frac{\langle u, x \rangle - \alpha}{\lVert u \rVert^2}u
\end{align}
and \eqref{eq_descent_property} turns into an equation, see, e.g., \cite{schopfer2008metric, schuster2012regularization}. 


%
%

By the following theorem we want to provide some tools that will later be essential to define the sequential subspace optimization techniques we use to obtain faster reconstructions of the stored energy function. 
Essentially, these techniques consist of sequential metric projections onto (intersections of) hyperplanes or stripes. By definition \ref{dfn:metrische Projektion} we already know that a metric projection onto a non-empty, closed convex set can be formulated as a minimization problem. 
The special case of metric projections onto intersections of hyperplanes is summarized in the following theorem. A proof  
can be found in \cite{schuster2012regularization} for the more general setting of Bregman projections in (convex and uniformly smooth) Banach spaces $X$ and $Y$.\\ 

\newpage

\begin{theorem}\label{satz:projmin}
\begin{enumerate}
\item[(a)] { Let $H(u_i, \alpha_i)$ be hyperplanes for $i=1,...,N$ with non-empty intersection 
\begin{align*}
H:= \bigcap_{i=1}^{N} H(u_i, \alpha_i).
\end{align*}
The projection of $x$ onto $H$ is given by 
\begin{align*}
P_H(x)=x-\sum \limits_{i=1}^{N} \tilde{t}_i u_i,
\end{align*}
where $\tilde{t}:=\left( \tilde{t}_1, ... , \tilde{t}_N \right) \in \mathbb{R}^N$ minimizes the convex function
\begin{align*}
h(t)=\frac{1}{2} \Big\lVert x - \sum \limits_{i=1}^N t_i u_i \Big\rVert^2 + \sum \limits_{i=1}^N t_i \alpha_i, \quad t=\left( t_1, ..., t_N \right) \in \mathbb{R}^N.
\end{align*}
The partial derivatives of the function $h(t)$ are given by
\begin{align}\label{eq_partial_derivatives_h}
\frac{\partial}{\partial t_j}h(t) = - \Big\langle u_j, x- \sum \limits_{i=1}^N t_i u_i \Big\rangle + \alpha_j.
\end{align}
If the vectors $u_i$, $i=1,\ldots ,N$, are linearly independent, $h$ is strictly convex and $\tilde{t}$ is unique. } 
\item[(b)] { Let $H_i:= H_\leq(u_i, \alpha_i)$, $i=1,2$, be two half-spaces with linear independent vectors $u_1$ and $u_2$. Then $\tilde{x}$ is the projection of $x$ onto $H_1 \cap H_2$ if $\tilde{x}$ satisfies the Karush-Kuhn-Tucker conditions for 
\begin{displaymath}
\min \limits_{z \in H_1 \cap H_2} \lVert z-x \rVert^2.
\end{displaymath}
The Karush-Kuhn-Tucker conditions are given by
\begin{align*}
\tilde{x}&=x-t_1 u_1-t_2 u_2 &&\textrm{for any } t_1, t_2 \geq 0, \\
\alpha_i &\geq \langle u_i, \tilde{x} \rangle,  &&i=1,2,\\
0 & \geq t_i \left( \alpha_i- \langle u_i,  \tilde{x} \rangle \right),  &&i=1,2.
\end{align*} }
\item[(c)]{ For $x\in H_>(u,\alpha)$ the projection of $x$ onto $H_\leq(u,\alpha)$ is given by  
\begin{align*}
P_{H_\leq(u,\alpha)}(x)=P_{H(u,\alpha)}(x)=x-t_+u
\end{align*}
with
\begin{align*}
t_+=\frac{\langle u, x \rangle - \alpha}{\lVert u \rVert^2} >0.
\end{align*} }
\item[(d)] { The projection of $x\in X$ onto the stripe $H(u,\alpha, \xi)$ is given by
\begin{align*}
P_{H(u,\alpha,\xi)}(x)= \begin{cases}
P_{H_\leq(u,\alpha+\xi)}(x) &\textrm{if } x\in H_>(u,\alpha+\xi),\\
x, 							&\textrm{if } x\in H(u,\alpha,\xi), \\
P_{H_\geq(u,\alpha-\xi)}(x) &\textrm{if } x\in H_<(u,\alpha-\xi).
\end{cases}
\end{align*} }
\end{enumerate}
\end{theorem}

Part (a) of Theorem \ref{satz:projmin} allows us to use tools from optimization (see also, e.g., \cite{nocedal2006numerical}) to determine the parameters $t=(t_1, ..., t_N) $.
The fact that the minimization of the function $h(t)$ corresponds to the projection onto the intersection of the hyperplanes $H(u_i,\alpha_i)$ for $i=1,...,N$ can be seen by taking a look at the partial derivatives \eqref{eq_partial_derivatives_h} of $h(t)$. Let us assume that the parameters $\tilde{t}=\left( \tilde{t}_1, ... , \tilde{t}_N \right)$ represent the local minimum of the function $h(t)$. Then, 
$$
\frac{\partial}{\partial t_j}h(\tilde{t}) = - \Big\langle u_j, x- \sum \limits_{i=1}^N \tilde{t}_i u_i \Big\rangle + \alpha_j=0.
$$ 
Since by definition we have
$$
  P_H(x) = x- \sum \limits_{i=1}^N \tilde{t}_i u_i,   
$$
we obtain
$$
\big\langle u_j, P_H(x) \big\rangle=\alpha_j
$$
for all $j=1,...N$, which shows that $P_H(x)=x-\sum \limits_{i=1}^{N} \tilde{t}_i u_i$ is an element of each hyperplane $H(u_i, \alpha_i)$, $i=1,...,N$ and, as a direct consequence, we have
\[   P_H (x)\in H.   \]

\begin{remark}
 If $F$ is a linear operator and the given data $y^{\delta}$ are noisy with noise level $0 \leq \lVert y^{\delta} - y \rVert \leq \delta$, then the solution set $\mathcal{M}_{Fx=y}$ of the linear operator equation $Fx=y$ is contained in the stripes $H(u,\alpha,\xi)$, where
 \begin{align*}
  u &:= F^* w\\
  \alpha &:= \big\langle w, y^{\delta} \big\rangle\\
  \xi &:= \delta \lVert w \rVert
 \end{align*}
 with arbitrary $w \in Y$, since for each $x \in \mathcal{M}_{Fx=y}$ we have
 \begin{displaymath}
  \begin{split}
   \left\lvert \left\langle u,x \right\rangle - \alpha \right\rvert &= \left\lvert \left\langle F^* w,x \right\rangle - \big\langle w, y^{\delta} \big\rangle \right\rvert \\
   &= \left\lvert \left\langle w,Fx - y^{\delta} \right\rangle \right\rvert = \left\lvert \left\langle w,y - y^{\delta} \right\rangle \right\rvert \\
   &\leq \delta \lVert w \rVert = \xi.
  \end{split}
 \end{displaymath}
 This observation is the basis to derive an iteration of the form
 \begin{displaymath}
  x_{n+1}^{\delta} = P_{H_n^{\delta}}\big(x_n^{\delta}\big), \quad n \in \mathbb{N},
 \end{displaymath}
 where $H_n^{\delta} := \bigcap_{i \in I_n} H(u_n^{\delta}, \alpha_n^{\delta}, \xi_n^{\delta})$ is the intersection of stripes containing the solutions of $Fx = y$. For each solution $x$, a reasonable choice of the parameters that define the stripes yields the descent property
 \begin{displaymath}
  \left\lVert x - x_{n+1}^{\delta} \right\rVert^2 \leq \left\lVert x - x_{n}^{\delta} \right\rVert^2 - C \big\lVert Fx_{n}^{\delta} - y^{\delta} \big\rVert^2.
 \end{displaymath}
 This property is used to show convergence and regularization properties of the method, see \cite{schuster2012regularization}.
\end{remark}


\subsection{RESESOP for nonlinear problems}
We turn to the regularizing sequential subspace optimization (RESESOP) technique for nonlinear inverse problems 
\begin{align} \label{Fx=y}
 F(x)=y, \quad F: \mathcal{D}(F) \subset X \rightarrow Y.
\end{align}
in Hilbert spaces $X,Y$ and noisy data $y^{\delta}$ with known noise level $\delta > 0$. The respective SESOP method that is applicable to unperturbed data can easily be derived by setting $\delta = 0$, see also \cite{WaldSchuster}.

In order to adapt the methods for linear operators to the nonlinear case, we must ensure that we project sequentially onto subsets of $X$ that contain the solution set
\begin{displaymath}
 \mathcal{M}_{F(x)=y} := \left\lbrace x \in \mathcal{D}(F) : F(x) = y  \right\rbrace
\end{displaymath}
of the operator equation \eqref{Fx=y}. In contrast to linear problems, we have to take into account the local character of nonlinear operators, i.e., we have to incorporate information on the local nonlinear behaviour of the forward operator into the definition of the stripes onto which we project in each iteration. To do this appropriately, we need the following assumptions on the operator $F$.

\vspace*{1ex}

Let $F:\mathcal{D}(F)\subset X \to Y$ be continuous and Fr\'{e}chet differentiable in an open ball 
\begin{align*}
B_\rho(x_0):=\left\lbrace x\in X : \lVert x-x_0 \rVert <\rho \right\rbrace \subset \mathcal{D}(F)
\end{align*}
around the starting value $x_0 \in \mathcal{D}(F)$ with radius $\rho > 0$ and let the mapping
\begin{align*}
B_\rho(x_0) \ni x \mapsto F'(x)
\end{align*}
from $B_\rho(x_0)$ into the space $L(X,Y)$ of linear and continuous mappings be continuous. \\
We assume there exists a solution $x^+ \in X$ of \eqref{Fx=y} that satisfies $x^+ \in B_\rho(x_0)$. 
This ensures that we start the iteration close to a solution, which is a mandatory requirement for nonlinear problems.\\
Furthermore, we assume that the forward operator $F$ satisfies the tangential cone condition
\begin{align} \label{kegelbedingung}
 \left\lVert F(x) - F(\tilde{x}) - F'(x)(x-\tilde{x}) \right\rVert \leq c_{\mathrm{tc}} \left\lVert F(x) - F(\tilde{x}) \right\rVert
\end{align}
with a positive constant $$0<c_\mathrm{tc}<1$$ and the estimate (continuity of the Fr\'{e}chet derivative)
\begin{align*}
\left\lVert F'(x) \right\rVert < c_F
\end{align*}
with $c_F >0 $ for all $x, \tilde{x}\in B_\rho(x_0)$.\\
We also assume that the operator $F$ is weakly sequentially closed. That is, for a weakly convergent sequence $\lbrace x_n\rbrace_{n \in \mathbb{N} }$ with $x_n \rightharpoonup x$ and $F(x_n)\to y$ holds
\begin{align*}
x\in \mathcal{D}(F) \quad \mathrm{and} \quad F(x)=y.
\end{align*}

If all these properties are fulfilled, we can formulate the RESESOP method as proposed in \cite{WaldSchuster} and obtain a regularization technique.\\

\begin{remark}
 The goal of general SESOP methods is to use multiple search directions $u_{n,i}^{\delta}$, $i \in I_n$, $\lvert I_N \rvert < \infty$, in each step $n \in \mathbb{N}$ of the iteration in combination with a regulation of the step width. We have $\mathcal{M}_{F(x)=y} \subset H(u_{n,i}^{\delta},\alpha_{n,i}^{\delta},\xi_{n,i}^{\delta})$ if we set
 \begin{align*}
  u_{n,i}^{\delta} &:= F'(x_i^{\delta})^*w^{\delta}_{n,i}\\
  \alpha_{n,i}^{\delta} &:= \big\langle w^{\delta}_{n,i}, F(x_i^{\delta}) -y^{\delta} \big\rangle - \big\langle F'(x_i^{\delta})^*w^{\delta}_{n,i}, x^{\delta}_i \big\rangle\\
  \xi_{n,i}^{\delta} &:= \lVert w^{\delta}_{n,i} \rVert \left( c_{\mathrm{tc}} \big(\lVert R^{\delta}_i \rVert + \delta \big) + \delta \right),
 \end{align*}
 see also \cite{WaldSchuster}. \\
 These definitions show that each hyperplane is related to the properties of $F$ close to the respective iterate. In particular, the noise level $\delta$ and the constant $c_{\mathrm{tc}}$ from \eqref{kegelbedingung} determine the width of the stripe: the higher the noise level and the larger the opening angle of the cone, the larger we have to choose the width of the stripe.\\
\end{remark}

\begin{figure}[H]
\centering
\def\svgwidth{250pt}
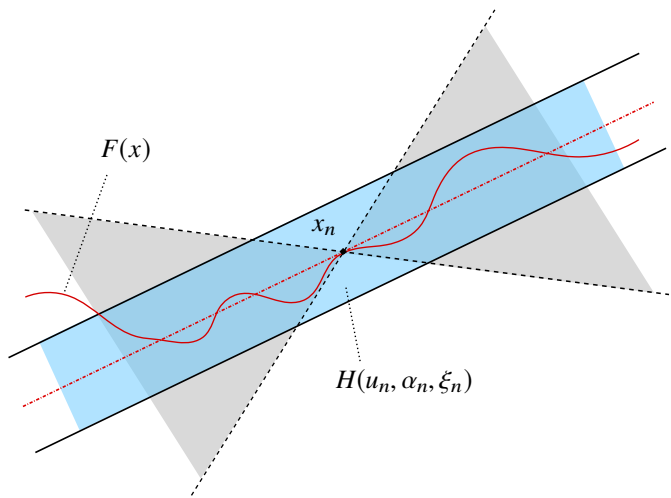
\caption{Illustration of a nonlinear function $F$ with stripe $H(u_n,\alpha_n, \xi_n)$}
\label{Abb: Streifenctc}
\end{figure}

Figure \ref{Abb: Streifenctc} illustrates the tangential cone condition \eqref{kegelbedingung} and its relevance for the choice of the stripes in the case $I_n := \lbrace n \rbrace$ for a function $F$ in two dimensions and exact data $y$. The graph of $F$ is plotted in red and for the point $x_n$ the linearization $F'(x_n)$ of $F$ in $x_n$ is represented by the red dotted line.  The graph is contained in the cone, determined by the tangential cone condition, highlighted in gray. The size of $c_{\mathrm{tc}}$ directly corresponds to the opening angle of the cone: The better $F$ is approximated by its linearization, the smaller is $c_{\mathrm{tc}}$ and thus also the opening angle of the grey cone. 
Figure \ref{Abb: Streifenctc} also shows that the cone condition can be used to define a stripe $H(u_n,\alpha_n, \xi_n)$ (marked in blue), such that the graph of $F$ is locally contained in $H(u_n,\alpha_n, \xi_n)$, i.e., in a neighborhood of $x_n$.\\

In the following we formulate the regularizing SESOP iteration for the special case of a single search direction per iteration, i.e., we set $I_n := \lbrace n \rbrace$ for all $n \in \mathbb{N}$. Furthermore, we define the $n$-th search direction as
\begin{displaymath}
 u_n^{\delta} := F'\big(x_n^{\delta}\big)^*\left(F\big(x_n^{\delta}\big) - y^{\delta}\right),
\end{displaymath}
such that we essentially obtain a Landweber-type method with an adaptation of the step size.
In comparison to the attenuated Landweber method, we thus have a dynamic relaxation parameter that adapts to the projection in each iteration step. Together with the discrepancy principle, we obtain a regularization method for which several convergence results could be shown (see \cite{WaldSchuster}). \\

\begin{algorithm} (RESESOP with one search direction) \\
We choose a starting value $x_0^\delta=x_0\in \mathcal{D}(F)$. For all $n \geq 0$ we select the search direction $u_n^\delta $ such that
\begin{align*} 
u_n^\delta&:= F'\big(x_n^\delta\big)^* w_n^\delta,\\
w_n^\delta&:= R_n^\delta:= F\big(x_n^\delta\big)-y^\delta.
\end{align*} 
We define the stripe $H_n^\delta$ by 
\begin{align*}
H_n^\delta:= H(u_n^\delta, \alpha_n^\delta, \xi_n^\delta)
\end{align*} 
with
\begin{align*}
\alpha_n^\delta&:= \langle u_n^\delta, x_n^\delta \rangle - \lVert R_n^\delta \rVert^2,\\
\xi_n^\delta&:= \lVert R_n^\delta \rVert \left( \delta + c_{tc} \left(\lVert R_n^\delta \rVert +\delta \right)   \right).
\end{align*}
As tolerance parameter for the discrepancy principle we choose 
\begin{align} \label{tauwahl}
\tau > \frac{1+c_\mathrm{tc}}{1-c_\mathrm{tc}}>1.
\end{align} 
As long as $\left\lVert R_n^\delta \right\rVert > \tau \delta$ is valid, we have
\begin{align} \label{xn_above_stripe}
x_n^\delta \in H_>(u_n^\delta, \alpha_n^\delta+\xi_n^\delta) 
\end{align}
and we calculate the new iterate $x_{n+1}^\delta$ by
\begin{align} 
x_{n+1}^\delta &:= P_{H(u_n^\delta, \alpha_n^\delta,\xi_n^\delta)}(x_n^\delta) = P_{H(u_n^\delta, \alpha_n^\delta+\xi_n^\delta)}(x_n^\delta) \label{it_proj}\\
 &=x_n^\delta - \frac{\langle u_n^\delta , x_n^\delta\rangle -\left(\alpha_n^\delta+\xi_n^\delta \right) }{\lVert u_n^\delta \rVert ^2} u_n^\delta.\label{it_expl}
\end{align}
\end{algorithm}

\begin{remark}
Note that due to \eqref{xn_above_stripe}, the iterate $x_n^\delta$ lies \emph{above} the stripe $H_n^{\delta}$ and, according to Theorem \ref{satz:projmin} (d), we obtain the identity \eqref{it_proj}. This projection is explicitly formulated in \eqref{it_expl}.\\
\end{remark}

The choice of $\tau$ in \eqref{tauwahl} depends strongly on the constant $c_\mathrm{tc}$ of the cone condition. The smaller $c_\mathrm{tc}$, the better the approximation of $F$ by its linearization. However, if $c_\mathrm{tc}$ is large, this also means that $\tau$ is large and the algorithm is usually stopped for larger residuals $\lVert R_n^{\delta} \rVert$. \\

For an analysis and a detailed discussion of general SESOP methods with multiple search directions in Hilbert and Banach space settings, we refer to the literature \cite{SchopferSchuster, schopfer2008metric, schuster2012regularization, WaldSchuster, Wald2018}.

%% file: nip_stripe.pdf_tex
\begingroup%
  \makeatletter%
  \providecommand\color[2][]{%
    \errmessage{(Inkscape) Color is used for the text in Inkscape, but the package 'color.sty' is not loaded}%
    \renewcommand\color[2][]{}%
  }%
  \providecommand\transparent[1]{%
    \errmessage{(Inkscape) Transparency is used (non-zero) for the text in Inkscape, but the package 'transparent.sty' is not loaded}%
    \renewcommand\transparent[1]{}%
  }%
  \providecommand\rotatebox[2]{#2}%
  \ifx\svgwidth\undefined%
    \setlength{\unitlength}{388.2bp}%
    \ifx\svgscale\undefined%
      \relax%
    \else%
      \setlength{\unitlength}{\unitlength * \real{\svgscale}}%
    \fi%
  \else%
    \setlength{\unitlength}{\svgwidth}%
  \fi%
  \global\let\svgwidth\undefined%
  \global\let\svgscale\undefined%
  \makeatother%
  \begin{picture}(1,0.73621844)%
    \put(0,0){\includegraphics[width=\unitlength]{nip_stripe.pdf}}%
    \put(0.13950054,0.54231677){\color[rgb]{0,0,0}\makebox(0,0)[lt]{\begin{minipage}{0.13635235\unitlength}\raggedright $F(x)$\end{minipage}}}%
    \put(0.4583378,0.4189515){\color[rgb]{0,0,0}\makebox(0,0)[lt]{\begin{minipage}{0.13427061\unitlength}\raggedright $x_n$\end{minipage}}}%
    \put(0.49543205,0.19142213){\color[rgb]{0,0,0}\makebox(0,0)[lt]{\begin{minipage}{0.35493248\unitlength}\raggedright $H(u_n,\alpha_n,\xi_n)$\end{minipage}}}%
  \end{picture}%
\endgroup%

%% file: num_res_arXiv.tex
\section{Numerical Results}
\label{sec:num_res}

In this section we present some numerical results to solve the inverse problem (IP) from Section \ref{sec:hyperelasticMaterials}. In all tests we use data that are simulated by solving the
initial boundary value problem \eqref{system:diff}--\eqref{system:diff3} using the $\theta$-method with respect to time and the Finite Element method in space. The resulting system of nonlinear equations is then solved by Newton's method. A detailed outline of the numerical forward solver for \eqref{system:diff} is contained in \cite{SEYDEL;SCHUSTER:17}.

The experimental setup for the numerical tests consists of a plate with measures $1\mathrm{m}\times 1\mathrm{m}$ and a thickness of $6.7\mathrm{mm}$. These measures can be numerically transferred to values of $\Omega = \left[-0.1 , 0.1\right]\times \left[-15,15\right]^2$. The plate is discretized using $5\times 31 \times 31$ knots with respect to $x$ and trilinear Finite Elements that are given by tensor products of linear B-splines. The time interval is given by $\left[0 \mu\mathrm{s}, 133 \mu \mathrm{s}\right]$, which we numerically represent as $\left[0,T\right]=\left[0,4\right]$. The time interval is discretized by $t_j = j\Delta t$, $j=0,\ldots,15$ and step size $\Delta t=0.25$. 
We assume that the plate is at rest at $t=0$ yielding $u_0=u_1=0$. The excitation signal $f(t,x)$ is chosen as a broad band signal that is emitted at the center of the plate acting in $x_3$-direction.
Again we refer to \cite{SEYDEL;SCHUSTER:17} for more details.

As already mentioned in Section \ref{sec:hyperelasticMaterials} the dictionary of stored energy functions $\{C_K: K=1,\ldots,N\}$ is defined as tensor products
\[  C_K (x,Y) = v_K (x) \hat{C}(Y) . \]
For our simulations we use the stored energy of a Neo-Hookean material model
\[  \hat{C} (Y) = c (I_1-3) + \frac{c}{\beta} (D^{-2\beta} - 1),   \]
where $I_1 = \|\nabla\varphi\|_F^2$, $D=\det (\nabla \varphi)$ and the constants are given by $\beta = \frac{3\nu - 2\mu}{6\mu}>0$ and $c=\frac{\mu}{2}>0$ with specific values
$\nu=68.6$ GPa and $\mu=26.32$ GPa taken from \cite{RAUTER;LAMMERING:15}. The functions $v_K$ are exactly the linear tensor product B-splines that are used for the Finite Element discretization
of the forward solver. Since linear tensor product B-splines have small compact support and represent a partition of unity, i.e.
\begin{equation}\label{part_unity}   \sum_{K=1}^N v_K(x) = 1,\qquad x\in\Omega,   \end{equation}
any defects can be appropriately modeled by coefficients $\alpha_K\not=1$ whereas for the undamaged plate we set $\alpha_K=1$, $K=1,\ldots,N$.

If we denote by $b_i$, $b_j$ the linear B-splines corresponding to the given discretizations in the $(x_2,x_3)$-plane, then we can simulate a delimation at the upper surface of the plate by defining the stored energy as
\begin{equation}\label{C_damage}  C(x,Y) := \sum_{i=0}^{30} \sum_{j=0}^{30} \alpha_{ij} b_i (x_2) b_j (x_3) \hat{C}(Y)\qquad \mbox{at } x_1=0.05   \end{equation}
and setting $\alpha_{ij}\not= 1$ for locations of the delamination.
Due to \eqref{part_unity}, $\alpha_{ij}=1$ corresponds to regions of the $(x_2,x_3)$-plane that are unaffected by the damage. Setting $\alpha_{ij}=1$ for all
$i$, $j$ yields $C(x,Y)=\hat{C}(Y)$ for all $x\in\Omega$ and thus models a homogeneous material. Note, that in \eqref{C_damage}
we use double indices in $\alpha_{ij}$ according to the tensor product structure of the Finite Elements $b_i\otimes b_j$, i.e., we have $\alpha_K=\alpha_{ij}$ with $K=31\cdot i + j$.

\begin{figure}[H]
\begin{tabular}{c c}
\includegraphics[height=4cm]{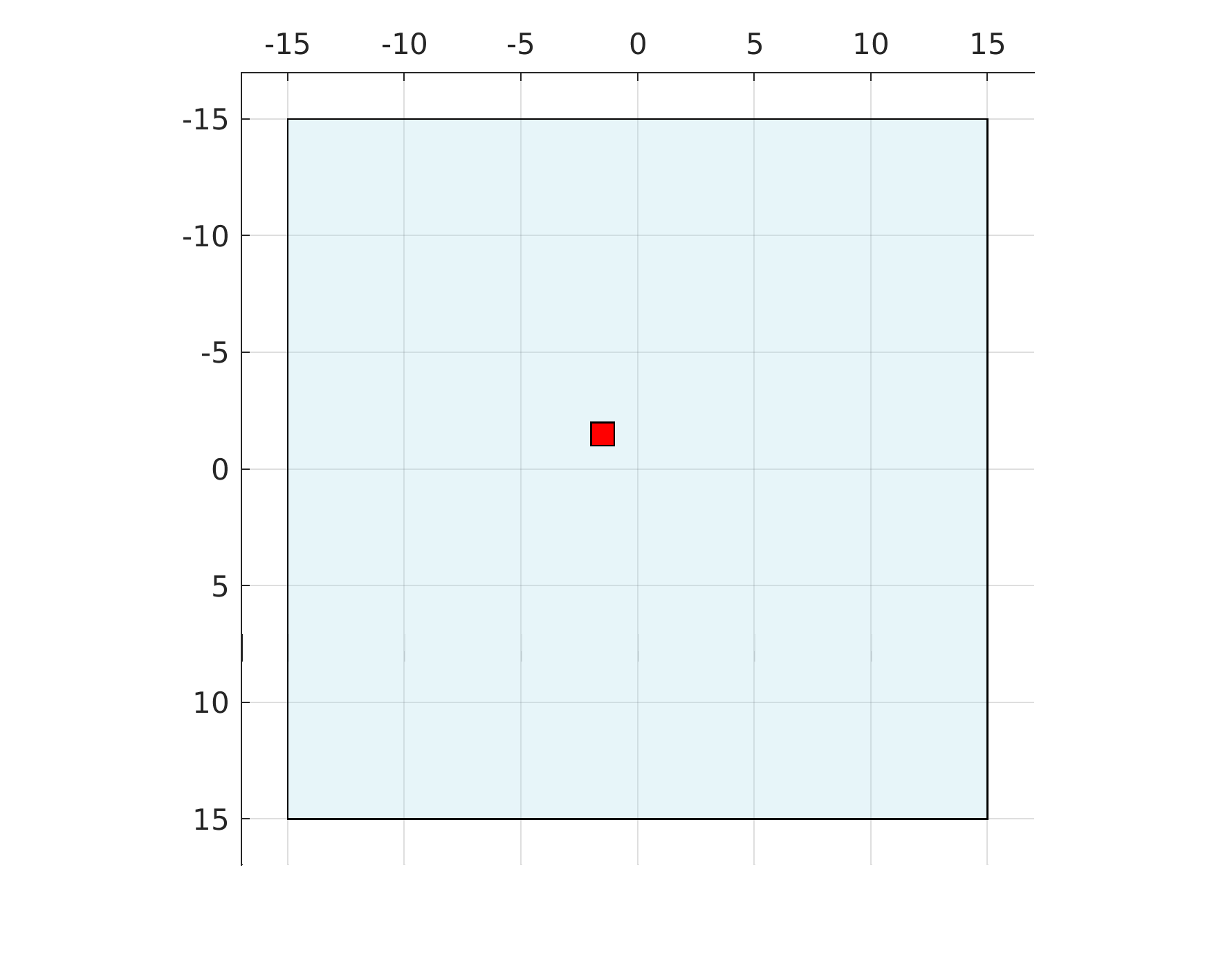}

&
\includegraphics[height=4cm]{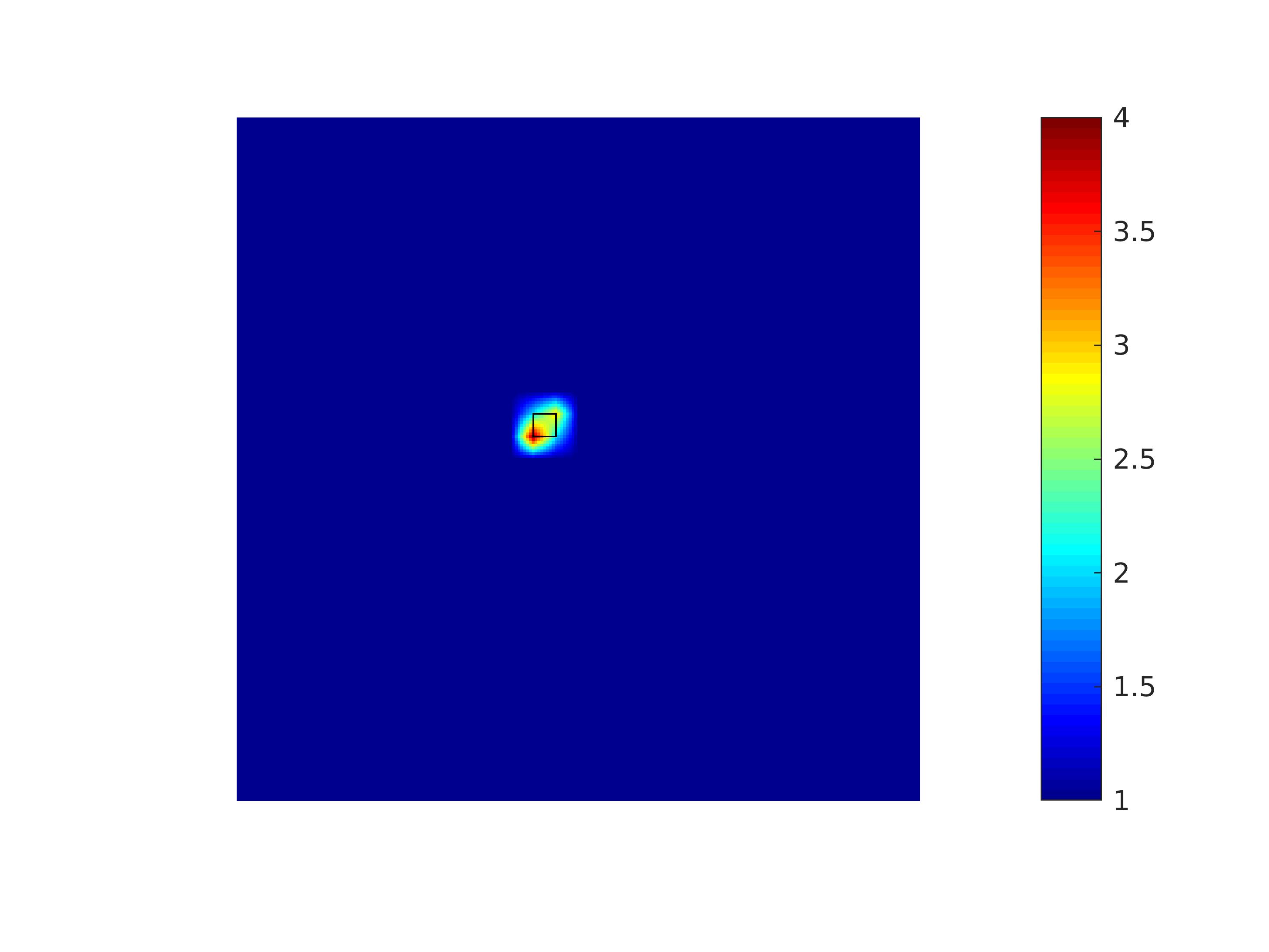}

\end{tabular}

\caption{Left picture: plate with damage at $\left(-1.5, -1.5 \right)$ (Experiment 1). Right picture: exact coefficient matrix $\alpha$ for experiment 1}
\label{abb:schadenV1_and_abb:schadenV1interpol}
\end{figure}

The first series of experiments examines a plate with a delamination whose center is located at $\left(x_2,x_3\right)=\left(-1.5, -1.5 \right)$, see Figure \ref{abb:schadenV1_and_abb:schadenV1interpol}. The corresponding coefficients $\alpha_{ij}$, $i,j\in \{0,\ldots,30\}$ in \eqref{C_damage} are given by
\begin{align*}
\alpha_{13,13}=2,\quad \alpha_{13,14}=3, \quad \alpha_{14,13}=4, \quad \alpha_{14,14}=2,
\end{align*}
and $\alpha_{i,j}=1$ elsewhere (Experiment 1). This setting for $\alpha_{ij}$ in fact corresponds to the damage in Figure \ref{abb:schadenV1_and_abb:schadenV1interpol} (left picture), which is emphasized in the right picture of Figure \ref{abb:schadenV1_and_abb:schadenV1interpol}  where the coefficient
matrix $\alpha = (\alpha_{i,j})_{i,j=0,\ldots 30}$ is plotted. There as well as in all reconstruction plots we apply linear interpolation to $\alpha$ to obtain a picture of higher resolution.
The inverse problem consists of computing the coefficient matrix $\alpha\in \mathbb{R}^{31\times 31}$ from full field data $u(t_j,x_m)$ where the discrete points $x_m$ correspond to the knots
of the Finite Element solver.\\

%
%

\begin{figure}[H]
\begin{tabular}{c c}
\includegraphics[height=4cm]{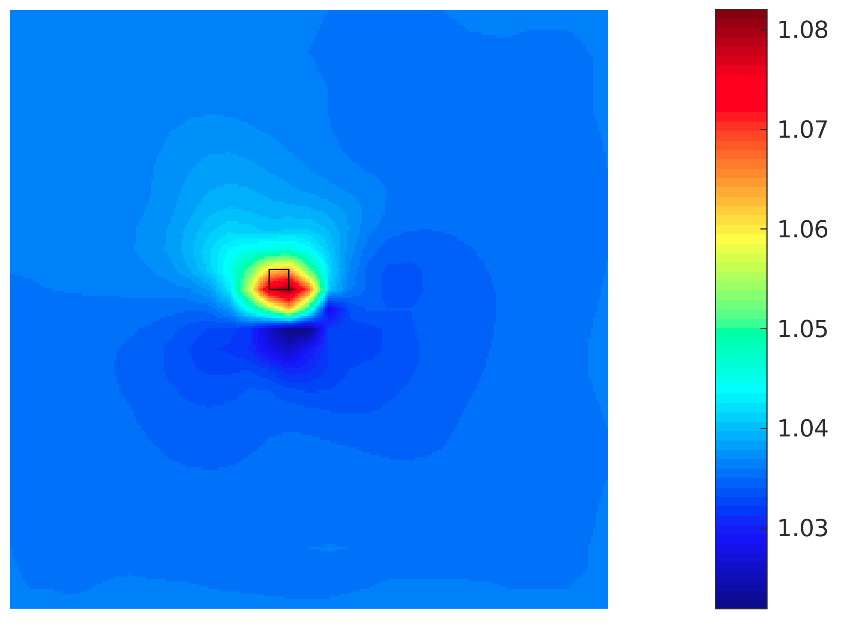}

&
\includegraphics[height=4cm]{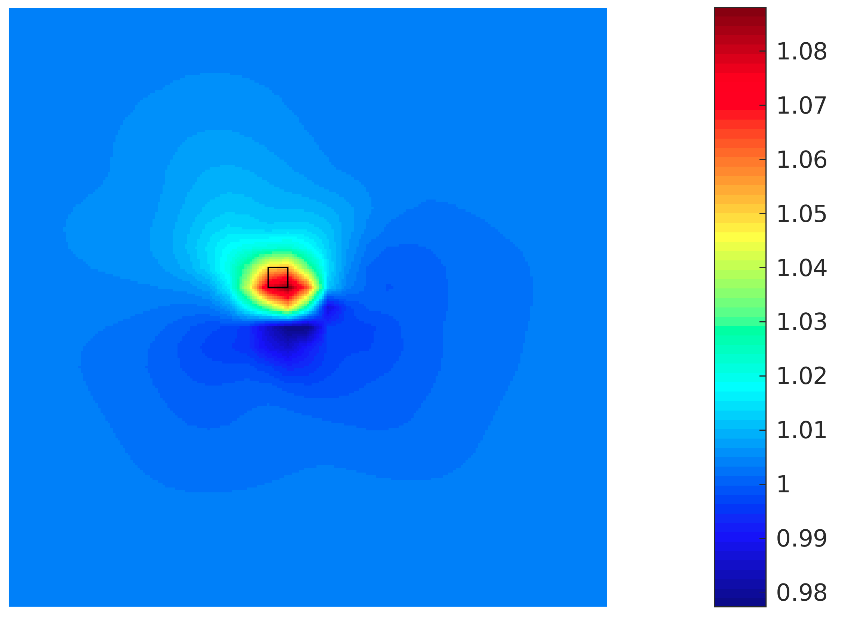}

\end{tabular}

\caption{Result of Experiment 1 after 200 iterations with the Landweber method (left) and after 9 iterations with the RESESOP method (right)}
\label{abb:v1sesopalpha_and_abb:v1landalpha}
\end{figure}

In the tests we compare different solution methods regarding the residual, the number of necessary iteration steps and computation time.
We implemented the Landweber iteration \eqref{eq_landweber} as well as RESESOP \eqref{it_expl} with the Landweber descent as single search direction and optimized step size in each iteration. Figure \ref{abb:v1sesopalpha_and_abb:v1landalpha} illustrates the results that are obtained after 9 iterations of RESESOP and 200 iterations of Landweber's method. 
The RESESOP iteration was stopped by the discrepancy principle, whereas the Landweber iteration was stopped before the discrepancy principle was fulfilled. In both cases the defect is detected at the correct location, but the coefficients $\alpha_{ij}$ are underestimated. We conclude that the same reconstruction quality is achieved with both methods but that RESESOP needs a significantly smaller number of iterations compared to the Landweber scheme. That means that RESESOP with only one search direction and optimized step size converges much faster than Landweber's method.

Next we compare the computing time that is needed for each iteration. One Landweber iteration needs 2.8 hours, resulting in a total computation time of 23 days until the discrepancy principle
is fulfilled. A RESESOP iteration takes 3 hours and thus a bit more than a Landweber step. But, since only 9 iterations are necessary to satisfy the discrepancy principle, the entire
reconstruction process only needs 27 hours in total. This means an acceleration by a factor of $\sim$51. We emphasize that (IP) is a high-dimensional parameter identification problem for a nonlinear hyperbolic system in time and space and thus belongs to the currently most challenging class of inverse problems at all.

\begin{figure}[H]
\centering
\includegraphics[height=6cm]{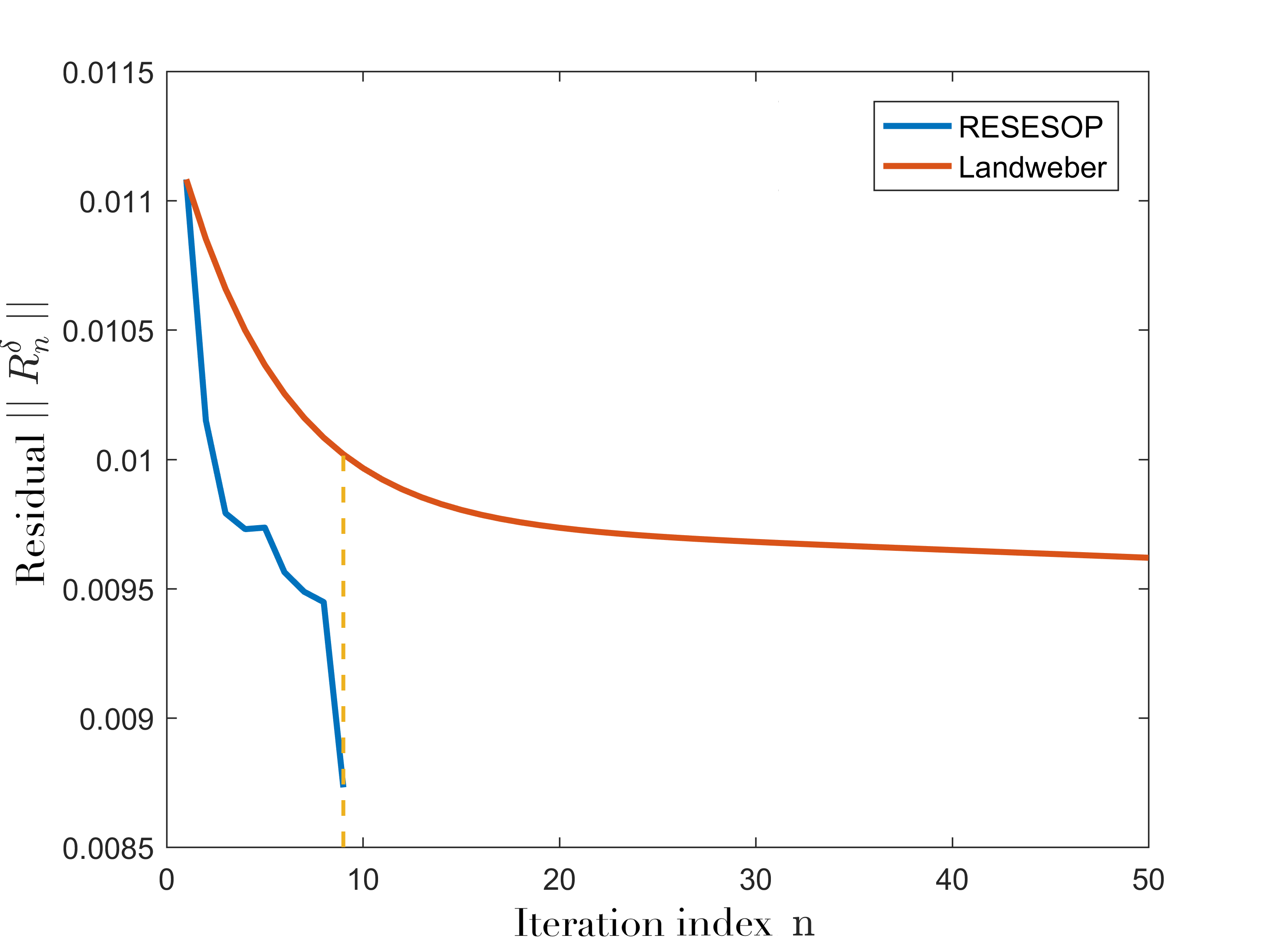}
\caption{Behavior of the respective residuals $ \lVert R_n^\delta \rVert$}
\label{abb:v1vergleichres}
\end{figure}

%
%
%

\begin{figure}[H]
\centering
\includegraphics[height=5cm]{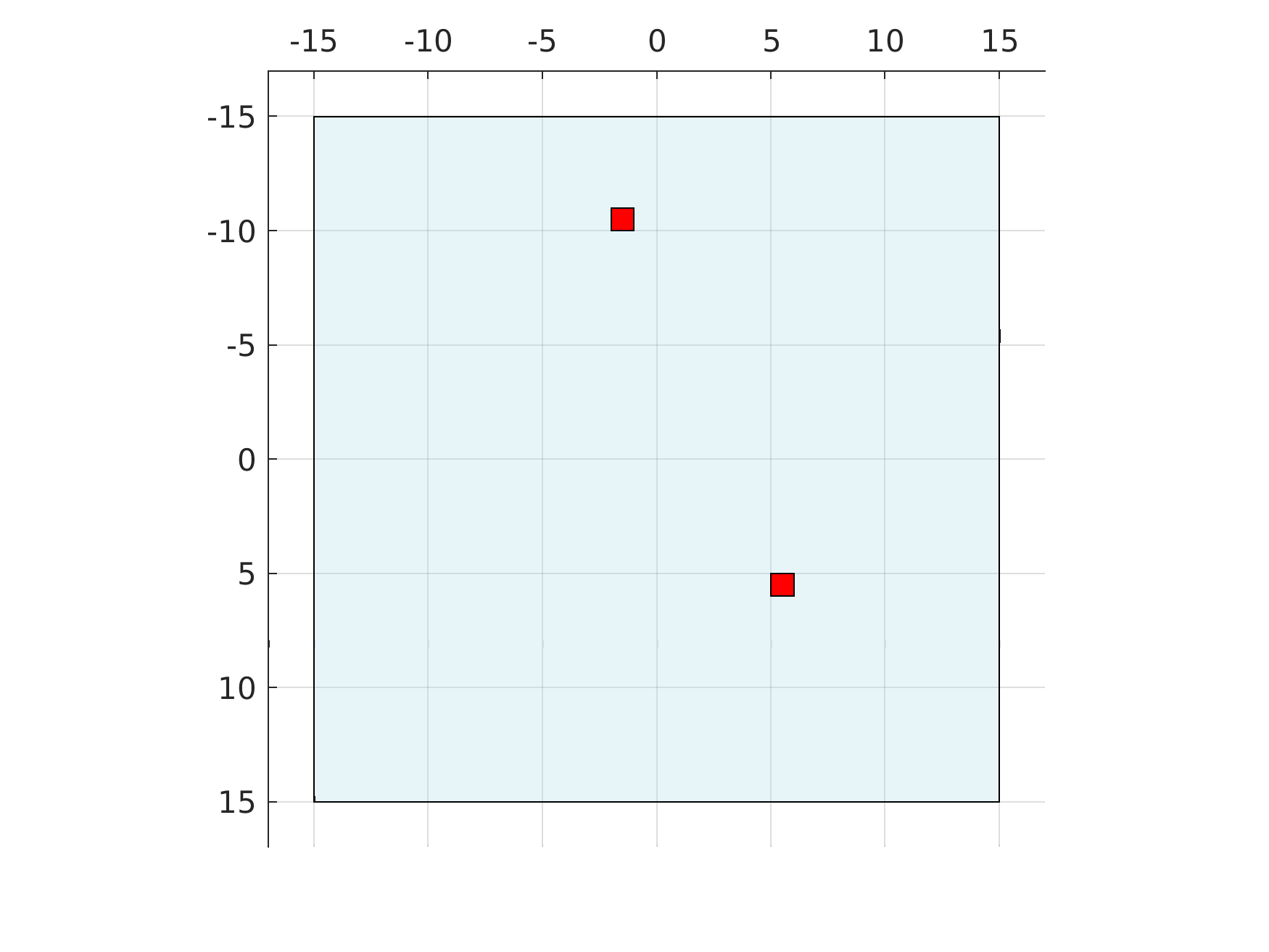}
\caption{Plate with damages A $\left( -1.5, -10.5\right)$ and B $\left(5.5, 5.5 \right) $ (Experiment 2) }
\label{abb:schadenV2}
\end{figure}

Figure \ref{abb:v1vergleichres} compares the residuals of the RESESOP technique and the Landweber method for Experiment 1. 
is satisfied. 
The red curve shows a typical behavior of the Landweber method. We see a strong decrease in the residual $\lVert R_n^\delta \rVert$ until iteration 15, followed by a very
slow decrease afterwards. This phenomenon is the reason why Landweber's method needs so much time until the discrepancy criterion is fulfilled.
We note also that the residual $\lVert R_n^\delta \rVert$ is not monotonically decreasing for RESESOP, in contrast to the Landweber iteration. The reason is that RESESOP
is constructed such that the sequence $\|x - x_n^\delta\|$ is monotonically decreasing, but not the sequence of residuals $\lVert R_n^\delta \rVert$, where $x$ denotes the exact solution of the underlying inverse problem and $x_n^\delta$ the $n$-th iterate for noisy data. This is in accordance with the analysis of the method outlined in \cite{WaldSchuster}.

In the second experiment we consider a setting consisting of two damages that are not located at the plate's center. Note that the center is also the region of wave excitation by $f(t,x)$. We assume that the damage which is closer to the center is the first to interact with the wave and thus is more pronounced in the reconstruction. 
The experimental setup is illustrated in Figure \ref{abb:schadenV2} (Experiment 2).\\

\begin{figure}[H]
\begin{tabular}{c c}
\includegraphics[height=4cm]{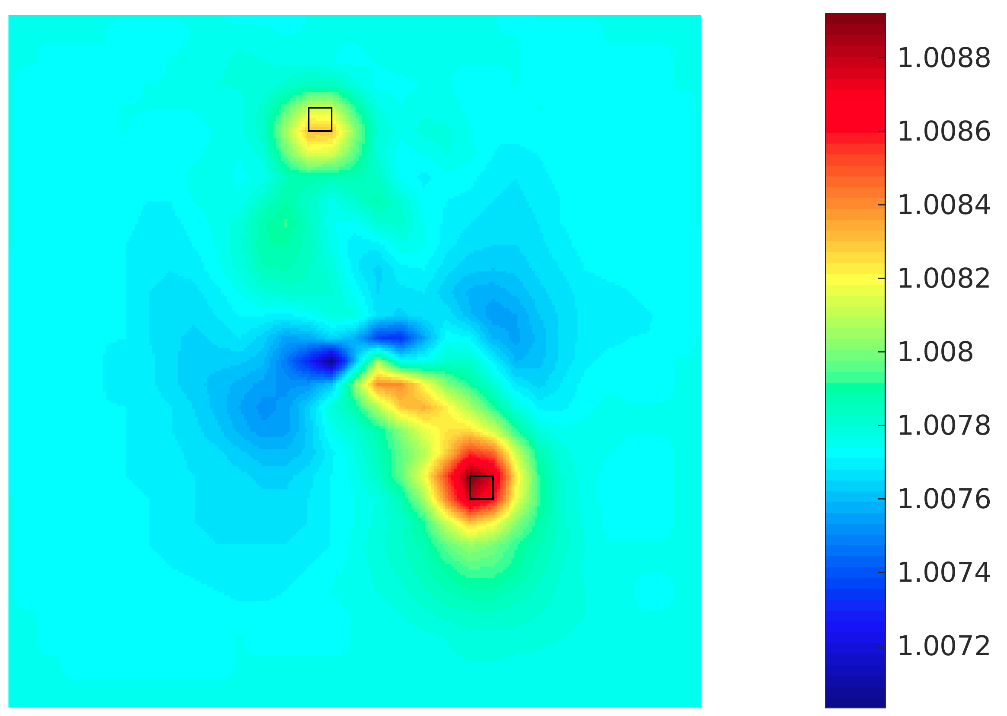}

&
\includegraphics[height=4cm]{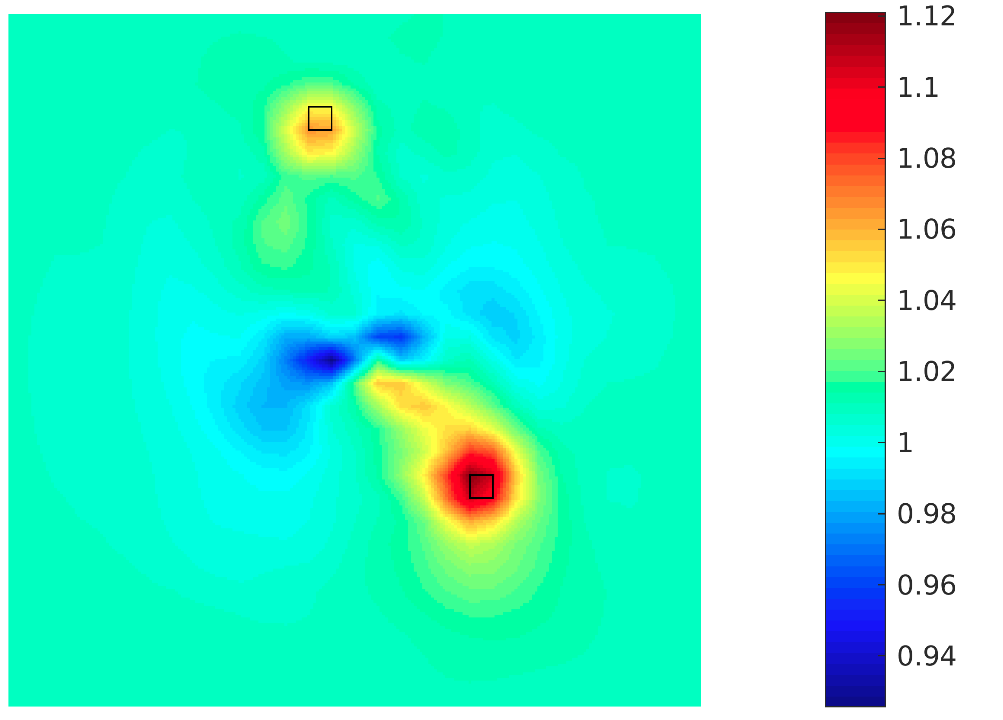}

\end{tabular}

\caption{Result of Experiment 2 after 50 iterations using the Landweber method (left) and after 17 iterations with the RESESOP method (right)}
\label{abb:v2landalpha_and_abb:v2sesopalpha}
\end{figure}

%
%
%

\begin{figure}[H]
\begin{tabular}{c c}
\includegraphics[height=4.5cm]{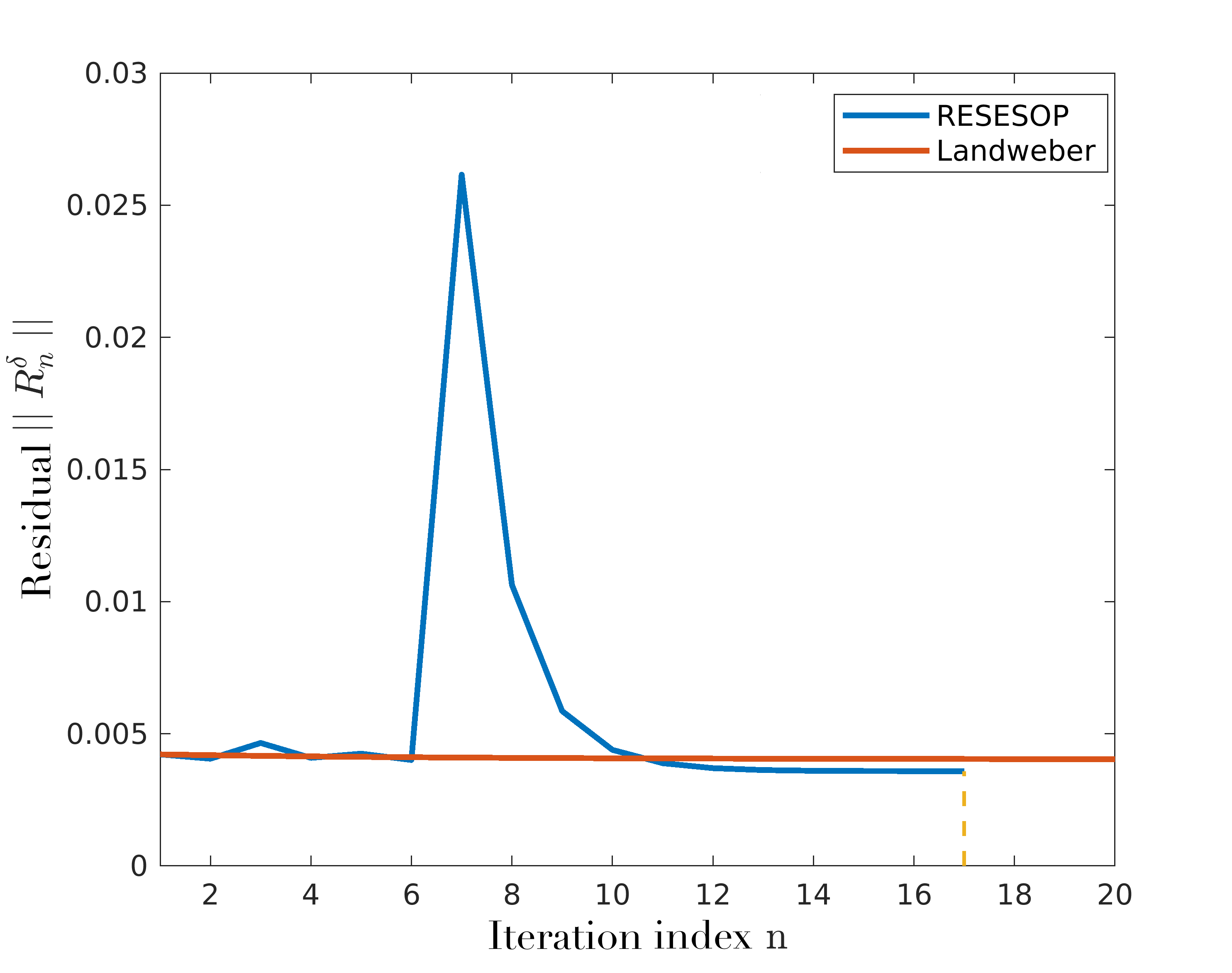}

&
\includegraphics[height=4.5cm]{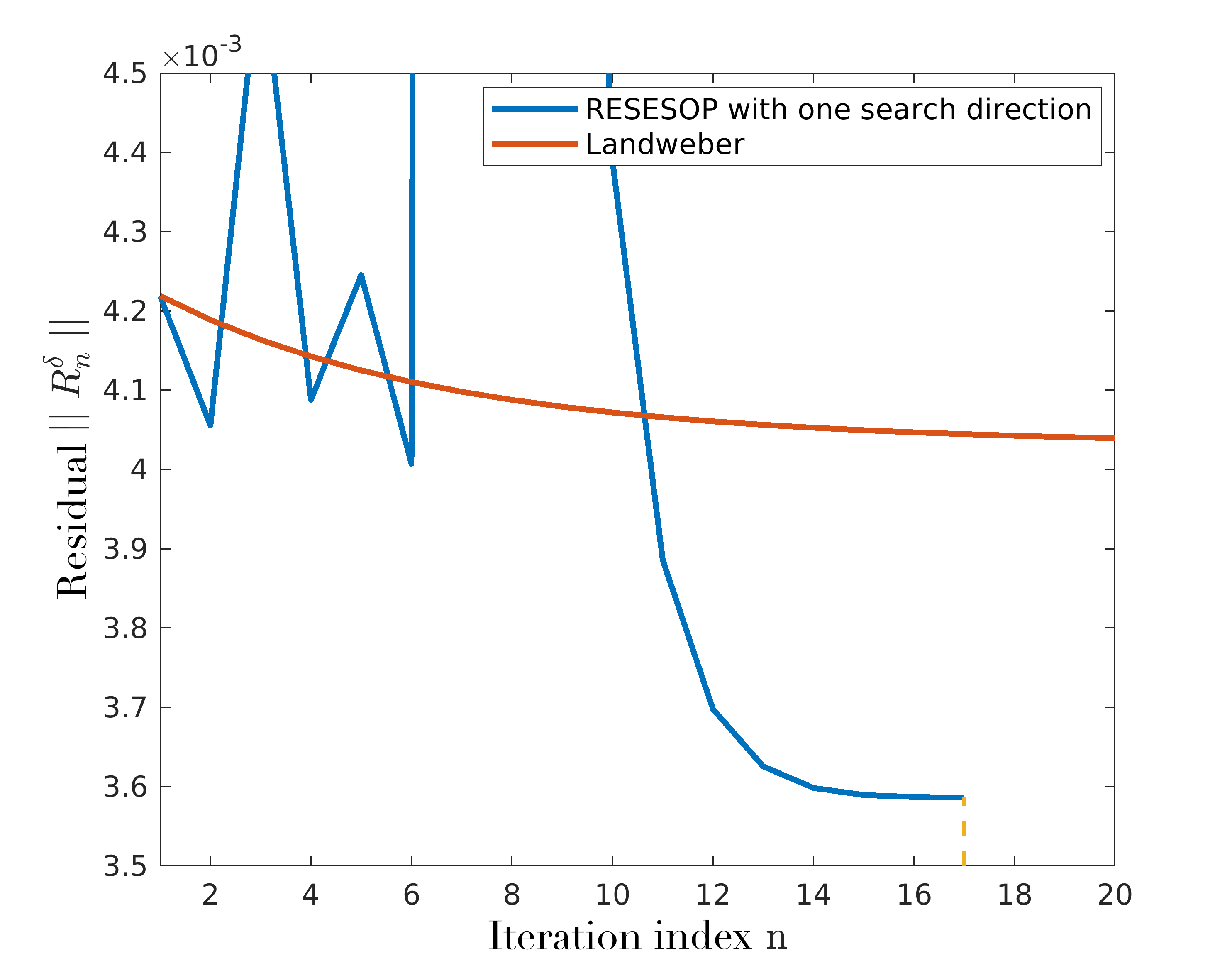}

\end{tabular}
\caption{Residuals $ \lVert R_n^\delta \rVert$ from the RESESOP and Landweber method from Experiment 2 (left picture) and a re-scaling of the $y$-axis (right picture)}
\label{abb:v2vergleichres_and_abb:v2rescut}
\end{figure}

Figure \ref{abb:v2landalpha_and_abb:v2sesopalpha} depicts the reconstruction from 50 iterations of the Landweber procedure (left picture). Then we terminated the iteration process because of its outrageous computation time. The values of the coefficient matrix $\alpha$ are contained in the very small interval $[1.0072, 1.0088]$. The situation is different
for the RESESOP technique. RESESOP stopped after iteration 17 according to the discrepancy principle. The result is visualized in Figure \ref{abb:v2landalpha_and_abb:v2sesopalpha} (right picture). The entries $\alpha_{i,j}$ of the coefficient matrix are contained in $[0.94, 1.12]$ making it easier to distinguish defects from undamaged parts of the structure. As expected the damage which is located closer to the center is more pronounced due to the excitation in the middle of the plate in both reconstructions.

In Figure \ref{abb:v2vergleichres_and_abb:v2rescut} we compare the residuals $\lVert R_n^\delta \rVert$ of the two methods applied to Experiment 2. 
The RESESOP iteration stops after iteration 17 according to the discrepancy principle, whereas the residual for the Landweber method seems to be almost constant. 
The right-hand plot in Figure \ref{abb:v2vergleichres_and_abb:v2rescut} shows a re-scaling to emphasize the oscillations of $\lVert R_n^\delta \rVert$ for RESESOP in the first few iterations as well as the monotonic decrease of $\lVert R_n^\delta \rVert$ for Landweber's method.
Furthermore both figures demonstrate again a faster convergence of RESESOP compared to the Landweber procedure.

\begin{figure}[H]
\centering
\includegraphics[height=8cm]{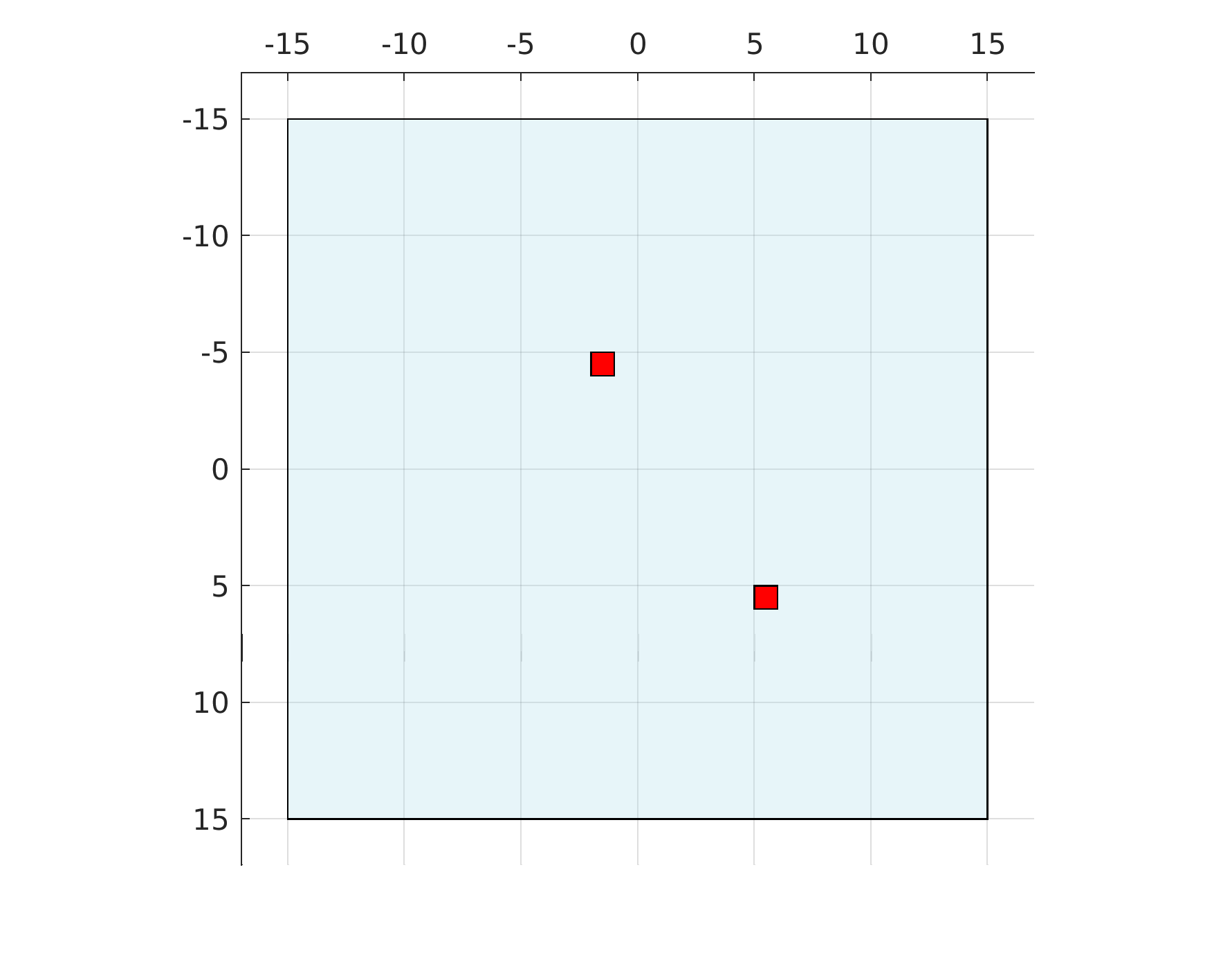}
\caption{Plate with damages A $\left( -1.5, -4.5\right)$ and B $\left(5.5, 5.5 \right) $ (Experiment 3) }
\label{abb:schadenV3}
\end{figure}

We consider a further numerical experiment where damage A is moved closer to the center of the plate compared to Experiment 2 and damage B remains fixed (Experiment 3).
This scenario is illustrated in Figure \ref{abb:schadenV3}. The corresponding coefficient matrix $\alpha$ remains unchanged, only the locations of the entries $\alpha_{i,j}$ are
adjusted to the damages A and B.

\begin{figure}[H]
\begin{tabular}{c c}
\includegraphics[height=4cm]{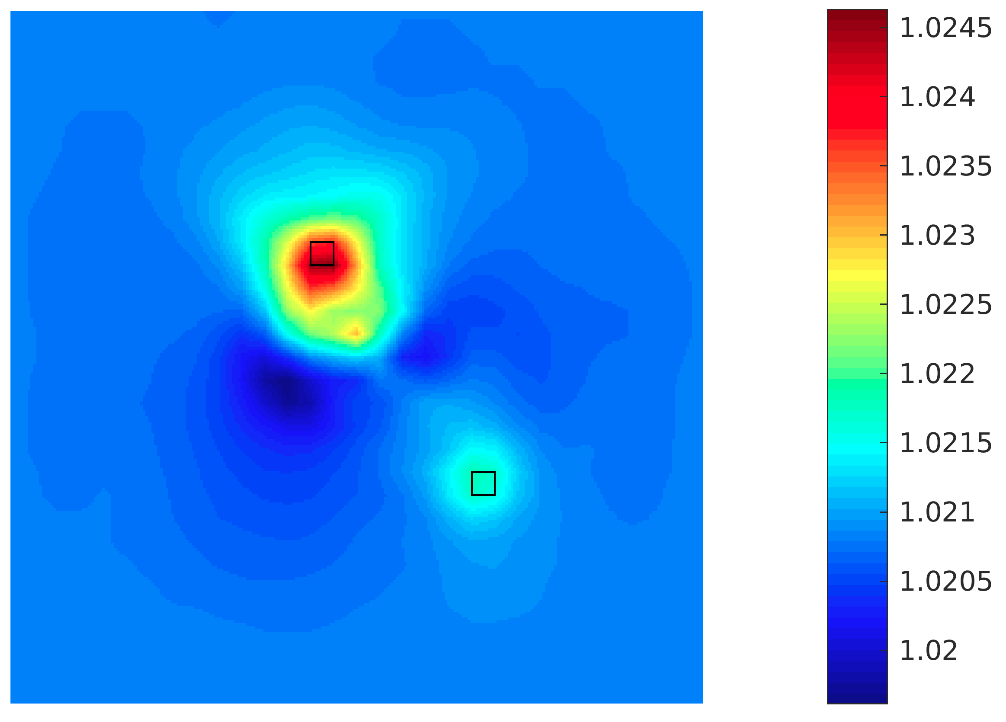}

&
\includegraphics[height=4cm]{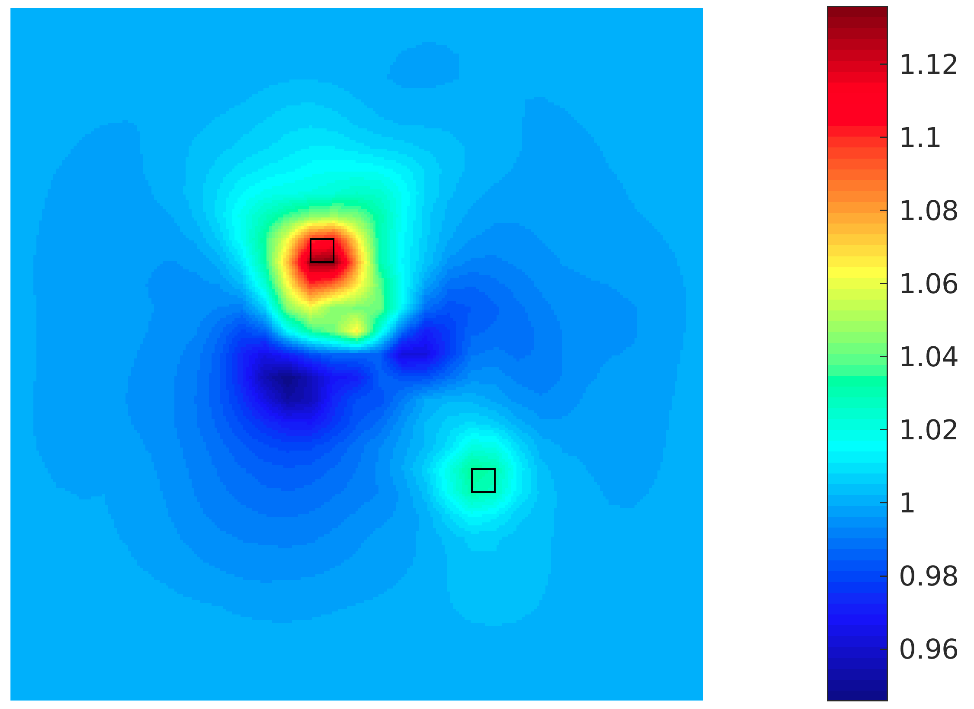}

\end{tabular}

\caption{Result of experiment 3 using the Landweber method after 50 iterations (left) and the RESESOP method after 14 iterations (right)}
\label{abb:v3alpha}
\end{figure}

Figure \ref{abb:v3alpha} shows the reconstructed coefficient matrix $\alpha$ using Landweber and RESESOP. In both cases the locations of the defects are accurately detected, while again the damage that is located closer to the center is highlighted stronger. The Landweber iteration has been stopped after 50 iterations (yielding 140 hours computation time) without having fulfilled the discrepancy principle. The RESESOP method, however, satisfied the discrepancy principle after 14 iterations (42 hours computation time) only, showing that it is significantly more efficient in spite of the additional computation time due to the step size optimization in each iteration.\\



The RESESOP technique outperforms the Landweber method in other respects as well. Considering the reconstructed values $\alpha_{i,j}$, we observe that the contrast in the Landweber reconstructions is very low, whereas an application of RESESOP results in larger differences of the absolute values.

\begin{figure}[H]
\begin{tabular}{c c}
\includegraphics[height=4.5cm]{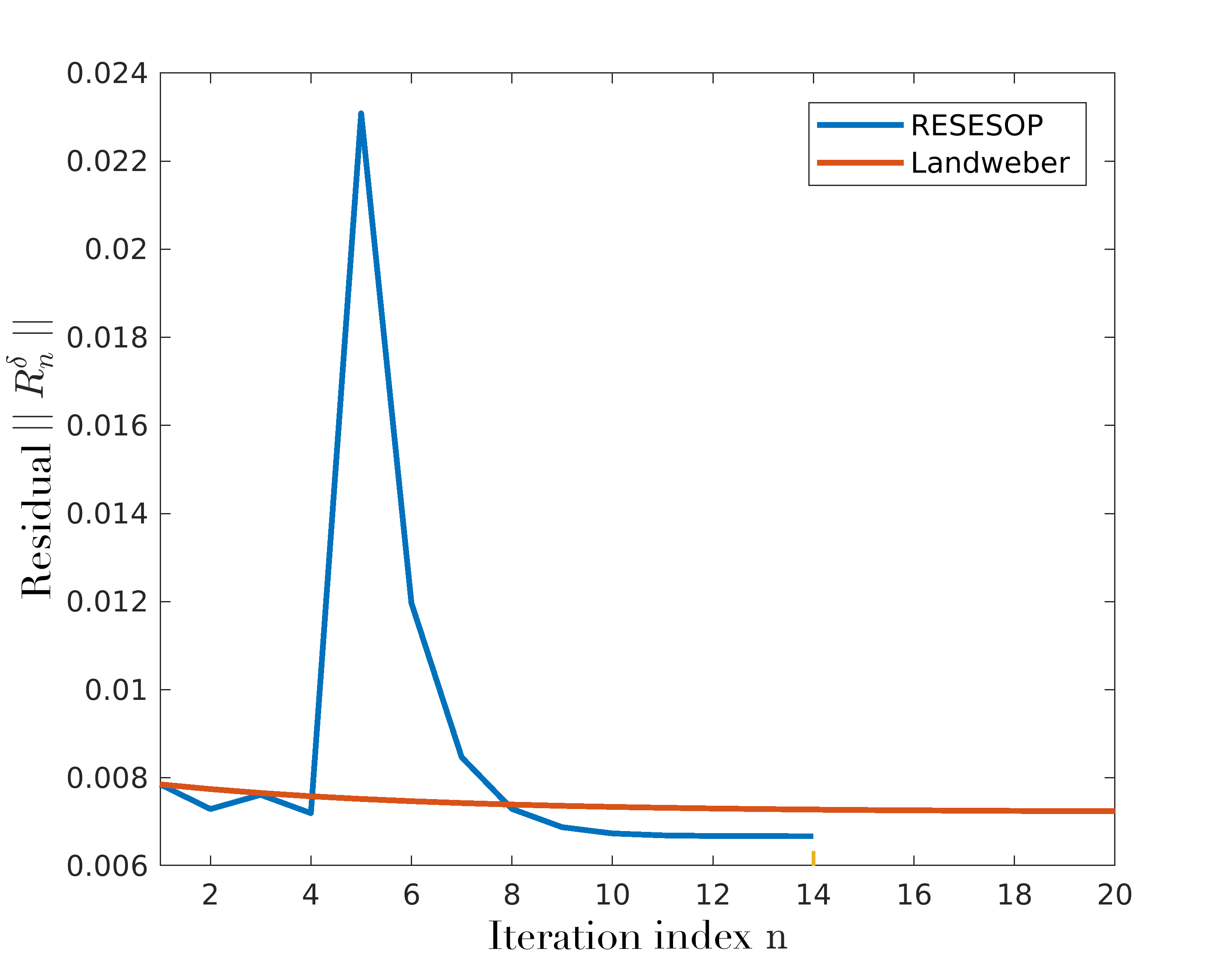}

&
\includegraphics[height=4.5cm]{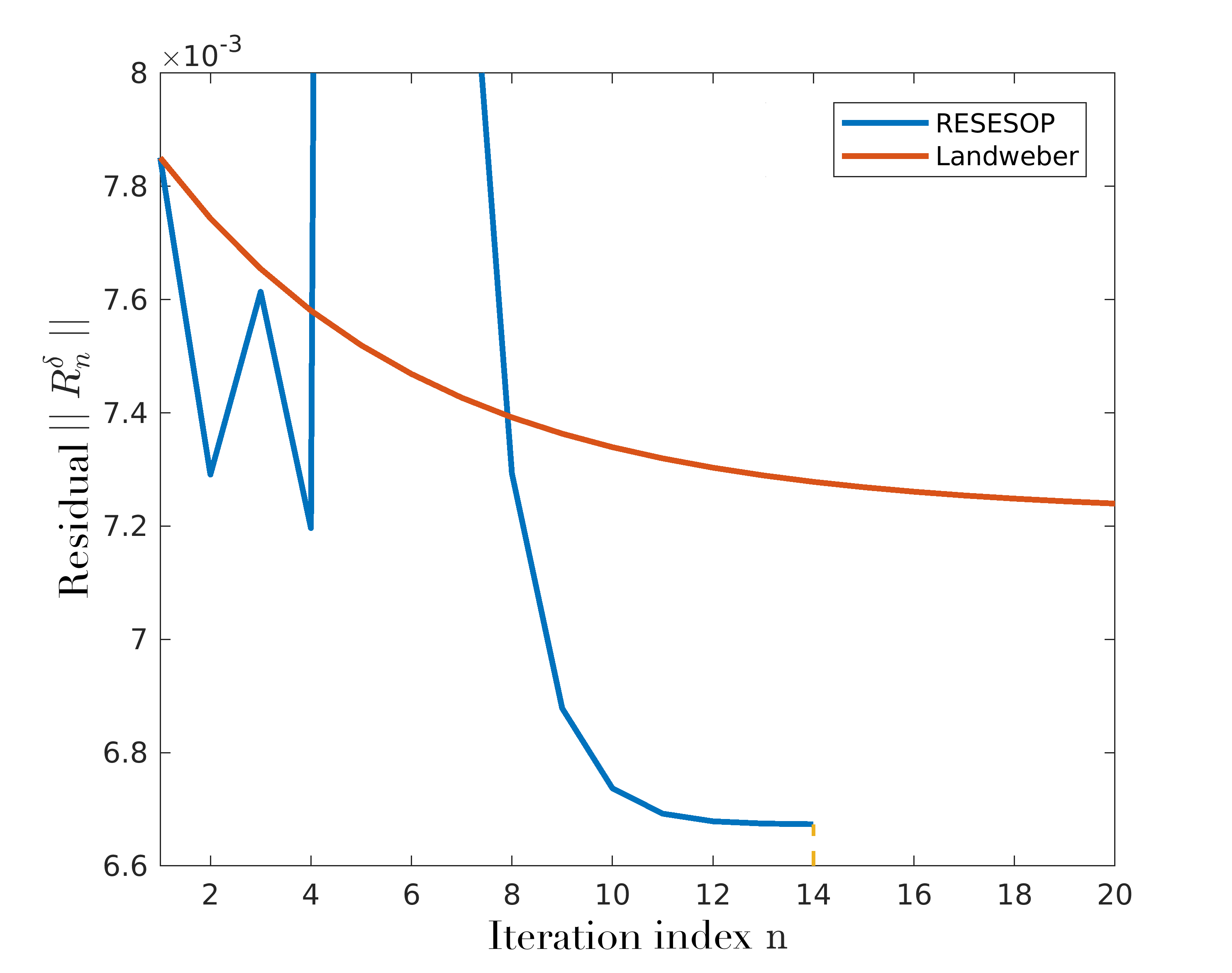}

\end{tabular}

\caption{Residuals in experiment 3 of the Landweber and the RESESOP method (left) and adapted scaling of the y-axis (right)}
\label{abb:v3res}
\end{figure}

Figures \ref{abb:v3res} and \ref{abb:v3alpha} show the residuals of the two methods when applied to Experiment 3. Similarly to Experiment 2, the figures clearly demonstrate the superiority of RESESOP compared to Landweber.\\

%% file: conclusion_arXiv.tex
\section{Conclusion}

We presented the performance of two different iterative regularization methods when applied to a high-dimensional inverse problem from the class of parameter identification
problems that is based on a system of nonlinear, hyperbolic differential equations equipped with initial and boundary values. The system describes
the propagation of elastic waves in a three-dimensional structure whose constitutive law is appropriately represented by a hyperelastic material model, i.e.
where the first Piola-Kirchhoff stress tensor is given as the derivative of the stored energy with respect to strain.
The nonlinearity allows also for large deformations. The considered inverse problem is the computation of the stored energy from measurements of the full displacement
field depending on space and time. Since the stored strain energy encodes virtually all essential mechanical properties of the structure on a macro-scale, it might yield
useful pointers for possible damages and thus might be important for simulations in the area of Structural Health Monitoring (SHM).

To solve this inverse problem we implemented the well-known Landweber method and Regularized Sequential Subspace Optimization (RESESOP) technique. The latter consists
of iterative metric projections onto hyperplanes that are determined by the used search directions, the nonlinearity of the forward mapping (via the constant in the tangential cone condition)
and the noise level. RESESOP uses in each iteration step a finite number of search directions where the Landweber direction, i.e. the negative gradient of the current
residual, is included. Using only one search direction, RESESOP coincides with Landweber where the step size is optimized to minimize the norm-distance of the current
iterate to a (locally unique) exact solution. Both numerical methods have been evaluated by means of three different damage scenarios for a Neo-Hookean material model and the usage of simulated measurement data. In all three cases RESESOP outperformes Landweber with respect to a faster convergence, a significant decrease of computation time and higher constrasts.

Future research could include model reduction techniques or the application of methods from Machine Learning. Both concepts
could help to achieve a further significant improvement with respect to computation time that is necessary for an implementation of the method in real-world SHM scenarios.

%% file: referenc_arXiv.tex
%
%